\documentclass[twoside]{article}

\usepackage[accepted]{aistats2022}
\usepackage[round]{natbib}

\usepackage{amssymb, amsmath, amsthm, graphicx, bbm}

\usepackage[dvipsnames]{xcolor}

\usepackage{enumitem}
\setlist{nolistsep}

\usepackage{algorithm}
\usepackage{algpseudocode}
\algnewcommand\algorithmicinput{\textbf{INPUT:}}
\algnewcommand\INPUT{\item[\algorithmicinput]}
\algnewcommand\algorithmicoutput{\textbf{OUTPUT:}}
\algnewcommand\OUTPUT{\item[\algorithmicoutput]}

\usepackage{hyperref}[]
\usepackage{cleveref}

\newtheorem{theorem}{Theorem}
\newtheorem{lemma}[theorem]{Lemma}
\newtheorem{proposition}[theorem]{Proposition}

\newtheorem{definition}{Definition}
\newtheorem{remark}{Remark}
\newtheorem{assumption}{Assumption}

\allowdisplaybreaks

\begin{document}

%

%

\twocolumn[

\aistatstitle{Optimal Partition Recovery in General Graphs}

\aistatsauthor{Yi Yu \And Oscar Hernan Madrid Padilla \And  Alessandro Rinaldo}

\aistatsaddress{ University of Warwick \And  University of California, Los Angeles \And Carnegie Mellon University } ]

\begin{abstract}
We consider a graph-structured change point problem in which we observe a random vector with piece-wise constant but otherwise unknown mean and whose independent, sub-Gaussian coordinates correspond to the $n$ nodes of a fixed graph.  We are interested in the localisation task of recovering the partition of the nodes associated to the constancy regions of the mean vector or, equivalently, of estimating the cut separating the sub-graphs over which the mean remains constant.  Although graph-valued signals of this type have been previously studied in the literature for the different tasks of testing for the presence of an anomalous cluster and of estimating the mean vector, no localisation results are known outside the classical case of chain graphs.  When the partition $\mathcal{S}$ consists of only two elements, we characterise the difficulty of the localisation problem in terms of four key parameters: the maximal noise variance $\sigma^2$, the size $\Delta$ of the smaller element of the partition, the magnitude $\kappa$ of the difference in the signal values across contiguous elements of the partition and the sum of the effective resistance edge weights $|\partial_r(\mathcal{S})|$ of the corresponding cut -- a graph theoretic quantity quantifying the size of the partition boundary.  In particular, we demonstrate an information theoretical lower bound implying that, in the low signal-to-noise ratio regime $\kappa^2 \Delta \sigma^{-2}  |\partial_r(\mathcal{S})|^{-1} \lesssim 1$, no consistent estimator of the true partition exists.  On the other hand, when $\kappa^2 \Delta \sigma^{-2}  |\partial_r(\mathcal{S})|^{-1} \gtrsim \zeta_n \log\{r(|E|)\}$, with $r(|E|)$ being the sum of effective resistance weighted edges and $\zeta_n$ being any diverging sequence in $n$, we show that a polynomial-time, approximate $\ell_0$-penalised least squared estimator delivers a localisation error -- measured by the symmetric difference between the true and estimated partition -- of order $ \kappa^{-2} \sigma^2 |\partial_r(\mathcal{S})| \log\{r(|E|)\}$. Aside from the $\log\{r(|E|)\}$ term, this rate is minimax optimal.  Finally, we provide discussions on the localisation error for more general partitions of unknown sizes.
\end{abstract}

\section{INTRODUCTION}

General graph-type data are ubiquitous in application areas, including social networks \citep[e.g.][]{odeyomi2020time}, neuroscience \citep[e.g.][]{khan2021effective}, climatology \citep[e.g.][]{mina2021network}, finance \citep[e.g.][]{liu2021dynamic} , biology \citep[e.g.][]{raimondi2021novel}, epidemiology \citep[e.g.][]{di2021psychological}, to name but a few.

In this paper, we consider a general graph-structured change point problem in which  we observe a random vector in $\mathbb{R}^n$ with unknown, piece-wise constant mean and whose independent sub-Gaussian coordinates correspond to the nodes of a fixed and arbitrary graph.  We are concerned with the localisation task of recovering the constancy regions of the mean vector in a manner that conforms to the topology of the underlying graph.  This is motivated by the emerging of network-type data, where the graphs are not necessarily chain graphs or grid graphs.  One example can be in epidemiology, people of study form a network and their interactions form edges.  It would be vital to accurately identify an abnormal cluster of people.  

To study this problem, we make the structural assumption that the partition associated with the mean vector specifies a multicut of the underlying graph of small weight, where the weight of each edge is its effective resistance, defined in \eqref{eq-def-effective-weight}. Informally, this assumption implies that the size of partition boundary is small or that the multicut is sparse relative to the topology of the graph.

The idea of weighting edges by the effective resistance in order to express the complexity of piece-wise constant graph-valued signals was originally put forward by \cite{fan2018approximate} for the related but different task of estimating a piece-wise constant signal over a graph in the squared error loss. In particular, the authors argue that such edge weighting scheme adapts to a varying degree of connectivity and spatial heterogeneity of the underlying graph. As a result, effective resistance provides  a natural and effective quantification of the complexity of a piece-wise constant signal, more so than naively assigning a unit weight to each edge, which amounts to an  $\ell_0$ complexity. As our results reveal, the key insight of \cite{fan2018approximate} extends to the localisation task of recovering  the partition associated to the constancy regions of the mean vector, as the edge weighting by effective resistance is essential to provide nearly-optimal localisation guarantees over arbitrary graph topologies.

To the best of our knowledge, graph-structured change point localisation problems for general graphs have not been considered in the literature. Indeed, change point analysis for piece-wise constant signals  traditionally assumes a total ordering of the coordinates to represent temporal changes, which can be trivially expressed using a chain graph (1-dimensional grid graph), or a $d$-dimensional grid graph.  In contrast, when the signal is allowed to conform to an arbitrary  graph in the manner described above, it is possible to obtain more sophisticated change point settings exhibiting a high degree of spatial complexity. However, due to the lack of a natural ordering of the coordinates in graph-structured mean change point problems,  virtually all the existing methodologies for change point localisation, which are  specifically designed to work in temporal settings, are inapplicable. To overcome this issue, we propose and analyse the properties of a change point estimator based on the approximate weighted $\ell_0$-penalised least squares methodology of \cite{fan2018approximate} -- a polynomial-time procedure deploying the $\alpha$-expansion algorithm of \cite{boykov2001fast}, to compute graph cuts.  We make the following contributions:

$\bullet$ We characterise the difficulty of the partition recovery task in terms of four critical parameters, which are allowed to change with the size and topology of the graph: the size $\Delta$ of the smallest element of the partition $\mathcal{S}^*$  of the nodes induced by the piece-wise constant means, the sub-Gaussian variance factor of the noise $\sigma^2$, the smallest magnitude $\kappa$ of the difference in the signal values across contiguous elements of the partition and the sum of the effective resistance edge weights for the multicut in the graph corresponding to   $\mathcal{S}^*$. Specifically, we show that in the low signal-to-noise regime in which $\kappa^2 \Delta \sigma^{-2}  |\partial_r(\mathcal{S}^*)|^{-1} \lesssim 1$, no procedure is guaranteed to estimate the partition $\mathcal{S}^*$ with a localisation error smaller than the order of $n$.

$\bullet$ We then focus on the case in which the partition induced by the mean value is known to contain only two elements. This setting has been considered in testing the presence of an anomalous clusters \citep[e.g.][]{arias2008searching, pmlr-v31-sharpnack13a}.  We remark that, since we do not assume that the constancy regions of $\mu$ correspond to connected sub-graphs, even this simplified case allows for multiple clusters.
We show that when 
    \[
        \kappa^2 \Delta \sigma^{-2}  |\partial_r(\mathcal{S}^*)|^{-1} \gtrsim \zeta_n \log\{r(|E|)\},
    \]
    where $r(|E|)$ is the sum of effective resistance weighted edges and $\zeta_n$ is any diverging sequence in $n$, the polynomial-time Algorithm \ref{alg-main} delivers a localisation error upper bounded by 
    \[
        \kappa^{-2} \sigma^2 |\partial_r(\mathcal{S}^*)| \log\{r(|E|)\}.
    \]
    We further prove that, aside from the $\log\{r(|E|)\}$ term, this rate is minimax optimal.  

$\bullet$ For partitions containing an unknown number of elements, we provide a general procedure for localisation and discuss localisation rates under a much stronger signal-to-noise ration condition. 

$\bullet$ We illustrate the strength of our methodology in a variety of experiments.

We emphasise that the localisation task over general graph-structured signals turns out to be fundamentally different, in both its theoretical and computational aspects,  from the detection task of testing for the presence of an anomalous cluster.  See the discussions in Sections~\ref{sec-connections} and \ref{sec-fan}. 

\subsection{Problem Setup}
We formalise our model and the localisation task.

\begin{assumption}[Model] \label{assume-model}
Let $G = (V, E)$ be a fixed  connected graph with vertex set $V = \{1, \ldots, n\}$ and edge set $E \subseteq V \times V$.  We observe a random vector $Y=(Y_1,\ldots, Y_n)^{\top}$ whose coordinates correspond to the vertices of $G$ and satisfy, for each $i\in V$, 
	\begin{equation}\label{eq-model}
		Y_i = \mu^*_i + \varepsilon_i,
	\end{equation}
	where $\mu^* = (\mu_1^*, \ldots, \mu^*_n)^{\top}$ is a  piece-wise constant mean vector and the errors  
	$\{\varepsilon_i\}_{i \in V}$ are centred i.i.d.~sub-Gaussian random variables with sub-Gaussian parameter $\sigma > 0$. We let $\mathcal{S}^* = \{S^*_1, \ldots, S^*_{K^*}\}$ be the partition of $V$ supporting the constancy regions of $\mu^*$. In detail, for some  $(f^*_1,\ldots,f^*_{K^*})^{\top} \in \mathbb{R}^{K^*}$, it holds that (i) for each $k \in \{1, \ldots, K^*\}$, $\mu^*_i = f^*_k$  for all $i \in S^*_k$;  and (ii) for any $k,l \in \{ 1,\ldots, K^*\}$ with $k \neq l$, if $\partial(S_k^*, S_l^*) = \{(i, j) \in E: \, i \in S_k^*, \, j \in S_l^*\} \neq \emptyset$, then $f_k^* \neq f_l^*$.
The mean vector $\mu^*$ and its associated partition $\mathcal{S}^*$ are unknown.
\end{assumption}

We remark that the whole graph $G$ is required to be connected, but each element of the partition is not required to induce a connected sub-graph, as it is customary in the literature on detection in graph-valued signals; see, e.g.~\cite{arias2011detection}.  This is demonstrated numerically in Case 4 in \Cref{sec-numeric}. 

Our goal is to recover the partition $\mathcal{S}^*$ accurately.  Specifically, under the structural assumption detailed in \Cref{assume-snr}, we seek an estimator of the partition $\widehat{\mathcal{S}}$ such that the Hausdorff distance between $\widehat{\mathcal{S}}$ and $\mathcal{S}^*$ normalised by the node size $n$ vanishes, as the node size goes to infinity, i.e.
    \begin{align}\label{eq-consistent-def}
        d_{\mathrm{H}}(\widehat{\mathcal{S}}, \mathcal{S}^*) = \max\left\{d_{\mathrm{H}_1}(\widehat{\mathcal{S}}, \mathcal{S}^*), \, d_{\mathrm{H}_1}(\mathcal{S}^*, \widehat{\mathcal{S}})\right\} = o(n),
    \end{align}
	as $n \to \infty$, where for any two partitions $\mathcal{A}$ and $\mathcal{B}$ of~$V$ we set $d_{\mathrm{H}_1}(\mathcal{A}, \mathcal{B}) = \max_{A \in \mathcal{A}} \min_{B \in \mathcal{B}} |A \triangle B|$ with $A \triangle B = (A \cup B) \setminus (A \cap B)$.

\subsection{Notation}\label{sec:notation}
For any $\delta > 0$, let $\delta \mathbb{Z} = \{m\delta\}_{m \in \mathbb{Z}}$.  For any $x \in \mathbb{R}$ and any $\delta > 0$, let $x^{\delta}$ be a closest value to $x$ in $\delta \mathbb{Z}$.  For any $\delta > 0$ and any $\mu \in \mathbb{R}^n$, a $\delta\mathbb{Z}$-expansion of $\mu$ is any other vector $\widetilde{\mu} \in \mathbb{R}^n$ such that there exists a single value $c \in \delta\mathbb{Z}$ satisfying that  for every $i \in \{1, \ldots, n\}$, either $\widetilde{\mu}_i = \mu_i$ or $\widetilde{\mu}_i = c$.  
For any partition  $\mathcal{S} = \{S_1, \ldots, S_K\}$, with $K \geq 2$, let $\partial(\mathcal{S}) = \{(i, j) \in E: \, i \in S_k, j \in S_l, k \neq l\}$.  For any edge weighting $w: \, E \to \mathbb{R}_+$ and any partition $\mathcal{S} = \{S_1, \ldots, S_K\}$, $K \geq 2$, let
    \begin{equation}\label{eq-partial-w-s-def}
        |\partial_w(\mathcal{S})| = \sum_{k = 1}^{K-1} \sum_{l = k+1}^K \sum_{i \in S_k, j\in S_l} w(i, j) \mathbf{1}\{(i, j) \in E\}
    \end{equation}
    and $w(|E|) = \sum_{(i, j) \in E} w(i, j)$.  A spanning tree $\mathcal{T}$ of $G$ is a sub-graph that is a tree which includes all the vertices of $G$.  Let $r:\, E \to \mathbb{R}_+$ be the effective resistance edge weights, that is 
    \begin{equation}\label{eq-def-effective-weight}
        r(i, j) = \frac{\#\mbox{ spanning trees that include }(i, j)}{\#\mbox{ spanning trees}}.
    \end{equation}
    For any connected graph $G$ with $n$ nodes and any vector $v \in \mathbb{R}^n$, we say $\mathcal{S}_v = \{S_1, \ldots, S_K\}$ is the partition induced by $v$, if and only if $\mathcal{S}_v$ is the partition with the smallest size that $v$ has constant values in each $S_k$, $k \in \{1, \ldots, K\}$.

\section{INFORMATION THEORETIC BOUNDS}\label{sec-lower-bound}

It is of our primary interest to understand the fundamental limits of localising the true partition of a general graph.  We first characterise the hardness of the problem using the following model parameters.

\begin{definition}
    With the notation in \Cref{assume-model}, when $K^* \geq 2$, let $\Delta = \min_k |S_k^*|$ and $\kappa = \min_{\substack{k \neq l \\ \partial(S_k^*, S_l^*) \neq \emptyset}} |f_k^* - f_l^*|$ be the minimal piece size  and the minimal jump size.
\end{definition}
Thus, for each distribution $P$ specifying a change point model as described in \Cref{assume-model}, we can associate the quantities $\Delta$, $\kappa$ and 
$|\partial_r(\mathcal{S}^*)|$. For simplicity, in our notation we omit the dependence on $P$.

We will show that the problem of recovering the partition is completely characterised by: $\kappa$ the minimal jump size, $\Delta$ the minimal size of constant signal, $\sigma$ the fluctuation of the noise and $|\partial_r(\mathcal{S}^*)|$ the connectivity between the partition, where $r$ consists of the effective resistance edge weights.  
Below, we firstly show in \Cref{prop-2} that, if $\kappa^2 \Delta \sigma^{-2} \leq c|\partial_r(\mathcal{S}^*)|$, where $c > 0$ is an absolute constant, then from an information-theoretic point of view, no algorithm is guaranteed to provide a consistent partition estimator, in the sense of \eqref{eq-consistent-def}. 

\begin{proposition}\label{prop-2}
Let $\mathcal{P}_1$ be the collection of joint distributions that
	\[
	\mathcal{P}_1 = \left\{P: \, \Delta = \min \left\{\lfloor c \kappa^{-2}\sigma^2|\partial_r(\mathcal{S}^*)| \rfloor, \, \lfloor n/4\rfloor\right\}\right\}.
	\]
	It holds that $\inf_{\widehat{\mathcal{S}}} \sup_{P \in \mathcal{P}_1} \mathbb{E}_P\{d_{\mathrm{H}}(\widehat{\mathcal{S}}, \mathcal{S}^*)\} \geq c_1 n$, where the infimum is taken over all possible partitions of $V$ and $c_1 > 0$ is an absolute constant.
\end{proposition}

Next, we demonstrate that, provided that $\kappa^2 \Delta \sigma^{-2} \geq C|\partial_r(\mathcal{S}^*)|$, where $C > 0$ is an absolute constant, the quantity $\kappa^{-2} \sigma^2 |\partial_r(\mathcal{S}^*)|$ is a minimax lower bound on the localisation error.

\begin{proposition}\label{prop-1}
Let $\mathcal{P}_2$ be the collection of joint distributions that $\mathcal{P}_2 = \left\{P: \, \kappa^2 \Delta \sigma^{-2} > C |\partial_r(\mathcal{S}^*)| \right\}$,
    where  $C > 0$ is an absolute constant.  It holds that $    \inf_{\widehat{\mathcal{S}}} \sup_{P \in \mathcal{P}_2} \mathbb{E}_P(d_{\mathrm{H}}(\widehat{\mathcal{S}}, \mathcal{S}^*)) \geq C_1 \kappa^{-2} \sigma^2 |\partial_r(\mathcal{S}^*)|$,
    where the infimum is taken over all possible partitions of $V$ and $C_1 > 0$ is an absolute constant.
\end{proposition}

\subsection{Connections with Other Literature}\label{sec-connections}

While our results for general graphs appear to be new, the special case of the chain graph has been thoroughly studied. Indeed, the partition recovery problem in the chain graph corresponds to the localisation task in the  change point literature.  In this case, it has been shown \citep[e.g.][]{wang2020univariate, verzelen2020optimal, chan2013detection} that for $K^* = 2$ (i.e.~$|\partial_r(\mathcal{S}^*)| = 1$), if $\kappa^2\Delta \sigma^{-2} \lesssim \log(n)$, then there is no algorithm guaranteed to be consistent; and provided that $\kappa^2 \Delta \sigma^{-2} \gtrsim \zeta_n$, a minimax lower bound on the error is $ \kappa^{-2} \sigma^2$.  Ignoring the logarithmic factors, we see that the results we have derived for general graphs in Propositions~\ref{prop-2} and \ref{prop-1} include the known results on chain graphs with $K^* = 2$.  Another special case of a general graph is the $d$-dimensional lattice.  To the best of our knowledge, the only known result comes from \cite{padilla2021lattice}, which studied localisation errors of regions constructed from unions of rectangles.  Under some stronger model assumptions, \cite{padilla2021lattice} showed that the estimation error, up to a logarithmic factor, is of order $\kappa^{-2} \sigma^2$. 

As for general graphs, we are not aware of any results in terms of the partition recovery accuracy in the existing work, but there has been a line of work on estimating the whole signal over the graphs, a.k.a.~de-noising, \citep[e.g.][]{jung2018network, kuthe2020engineering, fan2018approximate, han2019set, sharpnack2012sparsistency, hallac2015network,padilla2018dfs}.  We see partition recovery and de-noising as two closely related but different topics.  For instance, with chain graphs, it is well-understood that the fused lasso estimator \citep{tibshirani2005sparsity} is optimal for de-noising but sub-optimal in localising change points \citep{lin2017sharp}.  \cite{fan2018approximate} showed that the minimax optimal rate for de-noising a graph-structured signal with piece-wise constant mean $\mu^*$ is of order $\sigma^2 |\partial_r(\mathcal{S}^*)| \log(n/ |\partial_r(\mathcal{S}^*)|)$ and is achieved by iterating \Cref{alg-general}.  \Cref{prop-1} suggests that the relationship between the minimax optimal rate for de-noising and localisation  may be informally stated as 
    \[
        \mbox{localization rate} \asymp \kappa^{-2} \times \mbox{de-noising rate},
    \]
    a connection that has been previously established in the change point literature.  More discussions with \cite{fan2018approximate} will be provided in \Cref{sec-fan} after we present all theoretical results.

Another relevant and closely-related problem is that of detecting an abnormal cluster in a general graph \citep[e.g.][]{arias2008searching, addario2010combinatorial,  arias2011detection, arias2011global, hall2010innovated}.  For illustration purposes, we use the results in \cite{addario2010combinatorial} to compare with our findings.  Translated into our notation, the results in \cite{addario2010combinatorial} establish that, from a testing perspective, the detection boundary is 
    \begin{equation}\label{eq-addario-berry-result}
        \kappa^2 \Delta \sigma^{-2} \asymp \log(\# \mbox{ candidate sub-graphs}).
    \end{equation}
    Note that, without any additional assumption on the graph or on the class of candidate clusters, there are $2^n$ possible sub-graphs and the result \eqref{eq-addario-berry-result} reads as 
    \begin{equation}\label{eq-snr-boundary-detecti}
        \kappa^2 \Delta \sigma^{-2} \asymp n,
    \end{equation}
    which appears to be much stronger than the condition we require in \Cref{prop-1}.  
    One would then argue that this conclusion seems to contradict the conventional wisdom that testing is easier than estimation.  However, as what we will show later, this is in fact not the case, as the type of distributions and graph topologies for which  \eqref{eq-snr-boundary-detecti} holds form a subset of $\mathcal{P}_1$ in \Cref{prop-2}.

\begin{proposition}\label{prop-3}
Consider the family of distributions
	\[
		\mathcal{P}_3 = \left\{P: \kappa^2 \Delta \sigma^{-2} = cn, \, E = \{(i, j), 1 \leq i < j \leq n\}\right\},
	\]	
	where $c > 0$ is an absolute constant. It holds that $\inf_{\widehat{\mathcal{S}}} \sup_{P \in \mathcal{P}_3} \mathbb{E}_P \{d_{\mathrm{H}}(\widehat{\mathcal{S}}, \mathcal{S}^*)\}  \geq c_1 n$, where the infimum is taken over all possible partitions of $V$ and $c_1 > 0$ is an absolute constant.
\end{proposition}

Note that in \Cref{prop-3}, the definition $\mathcal{P}_3$ only includes complete graphs.  It follows from \Cref{lem-com-resis-weight} that $|\partial_r(\mathcal{S})|$ can be as large as of order $n$ for complete graphs.  \Cref{prop-3}, therefore, can be seen as a special case of \Cref{prop-2}.  

\begin{lemma}\label{lem-com-resis-weight}
For any complete graph $G$ with $n$ nodes and any nontrivial partition $\mathcal{S} = \{S_1, S_2\}$, $S_1, S_2 \neq \emptyset$, $|S_1| = n_1$ and $|S_2| = n_2$, it holds that $\partial_r(\mathcal{S}) = 2n_1 n_2 n^{-1}$.
\end{lemma}

The requirements in \eqref{eq-snr-boundary-detecti} and in $\mathcal{P}_3$ are rather pessimistic.  Since $\Delta \leq n$ by definition, they in fact assume that the signal strength $\kappa/\sigma$ should be at least of order $\sqrt{n}$, under which, off by a logarithmic factor, one may simply threshold the node values, completely ignore the graph structure and obtain a consistent detection and also partition.  This highlights the role of ``sparsity'' in the models: for detection problems, \eqref{eq-addario-berry-result} shows that a small class of candidate models requires a much smaller signal strength $\kappa/\sigma$; and for localisation problems, \Cref{prop-1} shows that a small cut $|\partial_r(\mathcal{S})|$ requires a much smaller signal strength $\kappa/\sigma$.

\section{CONSISTENT PARTITION RECOVERY}  \label{sec-upper-bound}
In this section, we introduce a polynomial-time algorithm for estimating the partition induced by the piece-wise constant mean vector as in \Cref{assume-model}. Towards that goal, we 
study the approximate solution to the least squared problem with weighted $\ell_0$ penalty considered by \cite{fan2018approximate}.   Given a graph $G = (V, E)$, observations $Y = (y_1, \ldots, y_n)^{\top} \in \mathbb{R}^n$, tuning parameters $\tau, \delta >0$ and an edge weighting $w \geq r$, consider the objective function
	\begin{equation}\label{eq-obj}
		F_w(\mu) = \frac{1}{2} \|Y - \mu\|^2 + \lambda \sum_{\{i, j\} \in E} w(i, j)\mathbf{1}\{\mu_i \neq \mu_j\}.
	\end{equation}
Throughout  $\widehat{\mu}$ denotes 
	a $(\tau, \delta \mathbb{Z})$-local-minimiser of \eqref{eq-obj}, defined next.

\begin{definition}\label{def1}
For $\delta > 0$ and $\tau \geq 0$, a $(\tau, \delta\mathbb{Z})$-local-minimiser of \eqref{eq-obj} is any $\mu \in \mathbb{R}^n$ such that for every $\delta\mathbb{Z}$-expansion $\widetilde{\mu}$ of $\mu$, $	F(\mu) - \tau \leq F (\widetilde{\mu})$.
\end{definition}
	
When the edge weights satisfy $w(i,j) = \mathbf{1}\{(i, j) \in E\}$, the function \eqref{eq-obj} is also known as the Potts functional \citep{potts1952some}.  In chain graphs, the minimizer of \eqref{eq-obj} has been succesfully used for change point detection and is thoroughly studied  \citep[e.g.][]{wang2020univariate, chan2013detection, verzelen2020optimal}.  An exact minimiser of \eqref{eq-obj} in chain graphs can be calculated in polynomial time \citep{friedrich2008complexity}, but this is not the case for general graphs.  \cite{fan2018approximate} proposed a polynomial-time algorithm, built upon the $\alpha$-expansion procedure \citep{boykov2001fast}. \cite{fan2018approximate}  also showed that for the de-noising task, a $(\tau, \delta\mathbb{Z})$-local-minimiser of \eqref{eq-obj} is optimal, when the weights are the effective resistance.

As our goal is to recover the partition, we consider a simpler  variant of Algorithm~1 in \cite{fan2018approximate}.  We first  focus on the case of  $K^*$ known and equal to $2$ -- a standard setting for detection \citep[see, .e.g,][]{arias2011detection} --  for which we obtain nearly optimal minimax rates. We then extend our results to the case of a general, unknown $K^*$. 

\subsection{The Case of Known \texorpdfstring{$K^* = 2$}{TEXT}}\label{sec-K=2}

The methodology we propose, detailed in \Cref{alg-main}, is a simple modification of Algorithm~1 in \cite{fan2018approximate}, which in turn is centred on the $\alpha$-expansion procedure \citep{boykov2001fast}.  We include it in \Cref{alg-sub} for completeness.

\begin{algorithm}[ht]
	\begin{algorithmic}
	    \INPUT $G = (V, E)$, $V = \{1, \ldots, n\}$, $\{Y_i, \, i \in V\}$,  $\delta$, $\tau$, $\lambda$,  $w$.
	    \State $\overline{Y} \leftarrow n^{-1} \sum_{i \in V} Y_i$, $\mathrm{FLAG} \leftarrow 0$.
		\State $(Y_{\min}, Y_{\max}) \leftarrow (\min_{i \in V} Y_i, \max_{i \in V}Y_i)$
		\State $(\widehat{\mu}, \mu_1) \leftarrow ((\overline{Y})^{\otimes n}, 0^{\otimes n}) \in \mathbb{R}^n \times \mathbb{R}^n$
		\For{each $c \in \delta\mathbb{Z} \cap [Y_{\min}, Y_{\max}]$}
		    \State $\widetilde{\mu} \leftarrow \alpha\mathrm{E}(\widehat{\mu}, Y, G, \lambda, w,c)$ \Comment{\Cref{alg-sub}}
			\If{$F_w(\widetilde{\mu}) \leq F_w(\widehat{\mu}) - \tau$}
				\State $\mathrm{FLAG} \leftarrow 1$.
				\If{$F_w(\widetilde{\mu}) \leq F_w(\mu_1)$}
					\State $\mu_1 \leftarrow \widetilde{\mu}$
				\EndIf
			\EndIf
		\EndFor
		\If{$\mathrm{FLAG} = 0$}
			\State $\mu_1 \leftarrow \widehat{\mu}$.
		\EndIf
		\OUTPUT $\mu_1$
		\caption{Variant of Algorithm~1 in \cite{fan2018approximate}.  $\mathrm{Potts}(G, Y, \delta, \tau, \lambda, w)$} \label{alg-main}
	\end{algorithmic}
\end{algorithm} 

\begin{algorithm}[ht]
    \begin{algorithmic}
        \INPUT $G = (V, E)$, $V = \{1, \ldots, n\}$, $Y, \mu \in \mathbb{R}^n$, $w: E \to \mathbb{R}_+$, $w \geq r$, $\lambda \geq 0$, $c$.
        \State Add two vertices $\texttt{s}$ and $\texttt{t}$
        \For{each $i \in V$}
            \State $\tilde{w}(s, i) \leftarrow (Y_i - c)^2/2$
            \State $\tilde{w}(t, i) \leftarrow (Y_i - \mu_i)^2 \mathbf{1}\{\mu_i \neq c\} + \infty\mathbf{1}\{\mu_i = c\}$
        \EndFor
        \For{each $(i, j) \in E$}
            \If{$\mu_i = \mu_j$}
                \State $\tilde{w}(i, j) \leftarrow \lambda w(i, j)\mathbf{1}\{\mu_i \neq c\}$
            \Else
                \State Add a new vertex $\texttt{a}_{i, j}$
                \State $E \leftarrow E \setminus \{(i, j)\} \cup \{(i, \texttt{a}_{i, j}), (j, \texttt{a}_{i, j}), (\texttt{t}, \texttt{a}_{i, j})\}$
                \State $\tilde{w}(i, \texttt{a}_{i, j}) \leftarrow \lambda w(i, j) \mathbf{1}\{\mu_i \neq c\}$
                \State $\tilde{w}(j, \texttt{a}_{i, j}) \leftarrow \lambda w(i, j) \mathbf{1}\{\mu_j \neq c\}$
                \State $\tilde{w}(\texttt{t}, \texttt{a}_{i, j}) \leftarrow \lambda w(i, j)$
            \EndIf
        \EndFor
        \State Find the minimum $\texttt{s}$-$\texttt{t}$-cut $(S, T)$ of the newly constructed graph based on weights $\tilde{w}$ such that $\texttt{s} \in S$ and $\texttt{t} \in T$.
        \For{each $i \in V$}
            \State $\widetilde{\mu}_i \leftarrow c \mathbf{1}\{i \in T\} + \mu_i \mathbf{1}\{i \in S\}$
        \EndFor
        \OUTPUT $\widetilde{\mu}$
        \caption{$\alpha$-expansion from \cite{boykov2001fast}. $\alpha\mathrm{E}(\mu, Y, G, \lambda, w,c)$. } \label{alg-sub}
    \end{algorithmic}
\end{algorithm}

In \cite{fan2018approximate}
\Cref{alg-sub} is called iteratively  until the objective function $F_w(\cdot)$ is not improved.  This results in a partition  of potentially more than two subsets of nodes, which is desirable for the purpose of signal estimation.  Since our goal is to recover a partition with two elements, 
our strategy is different.  In the language of \Cref{alg-main}, with the initialiser $\widehat{\mu}_i = \overline{Y}$ for all $i$, \Cref{alg-sub} is summoned repeatedly in order to find a value $c \in \delta\mathbb{Z}$ and a subset $V_1 \subset V$ such that $F_w(\widehat{\mu})$ reaches its minimum, where $\widehat{\mu}_i = c$, $i \in V_1$, and $\widehat{\mu}_i = \overline{Y}$, $i \in V \setminus V_1$.

To evaluate the computational cost of \Cref{alg-main}, we adapt the proof of Proposition~2.2 in \cite{fan2018approximate}.  Note that there are at most $(Y_{\max} - Y_{\min})/\delta$ different choices of $c$.  For each value $c$, the augmented graph in \Cref{alg-sub} has $O(|E|)$ vertices and edges.  Solving minimum s-t cut using either the Edmonds--Karp or Dinic algorithm requires time $O(|E|^3)$.  This implies that the computational cost of \Cref{alg-main} is of order $O\{|E|^3(Y_{\max} - Y_{\min})\delta^{-1}\}$.  We acknowledge that the computational cost is high, despite being polynomial-time.  We do not know whether there exist faster algorithms enjoying the same theoretical optimality.  In chain graphs, faster algorithms exist, yet theoretically sub-optimal.

Note that edge-weights are used in the $\alpha$-expansion subroutine and it means that  \Cref{alg-main} is edge-weights dependent.  It is naturally the case that the partition recovery performance of \Cref{alg-main} also depends on the choice of edge-weights.  

We state our result in generality by assuming that the edge weights are point-wise larger than the effective resistance, as done in \cite{fan2018approximate}. 

\begin{assumption} [Signal-to-noise ratio] \label{assume-snr}
For a certain choice of edge-weights $w:\, E \to \mathbb{R}_+$, with $w \geq r$, assume that there exists a constant $C_{\mathrm{SNR}} > 0$ such that $\kappa^2 \Delta \geq C_{\mathrm{SNR}} \sigma^2 |\partial_w(\mathcal{S}^*)| \log\{w(|E|)\} \zeta_n$, where $\zeta_n$ is any arbitrarily diverging sequence, as the number of node $n$ grows unbounded.
\end{assumption}

\begin{theorem}[$K^* = 2$ and $K^*$ is known]\label{thm-1}
Let $\delta \leq \sigma/\sqrt{n}$, $\tau \leq \sigma^2$ and $\lambda = C_{\lambda}\sigma^2\log\{w(|E|)\}$, where $C_{\lambda} > 0$ is an absolute constant.  Assume that $K^* = 2$ and $K^*$ is known.  Under Assumptions~\ref{assume-model} and \ref{assume-snr}, letting $\widehat{\mu}$ be an output of \Cref{alg-main},  with $w \geq r$, it holds with probability at least $1 - \{w(|E|)\}^{-c}$ that
    \begin{equation}\label{eq-thm1-result}
		d_{\mathrm{H}}(\widehat{\mathcal{S}}, \mathcal{S}^*) \leq C \kappa^{-2} \sigma^2 |\partial_w(\mathcal{S}^*)| \log\{w(|E|)\},
	\end{equation}
	where $c, C > 0$ are absolute constants and $\widehat{\mathcal{S}}$ is the partition induced by $\widehat{\mu}$.
\end{theorem}

Recalling the consistency definition \eqref{eq-consistent-def} and in view of \Cref{assume-snr}, we see that the estimation error satisfies
    \begin{align*}
        n^{-1} d_{\mathrm{H}}(\widehat{\mathcal{S}}, \mathcal{S}^*) & \leq Cn^{-1} \kappa^{-2}\sigma^2 |\partial_w(\mathcal{S}^*)| \log\{w(|E|)\} \\
        & \leq C C_{\mathrm{SNR}} \zeta_n^{-1} n^{-1} \Delta \to 0,
    \end{align*}
    which implies that the output of \Cref{alg-main} provides a consistent partition recovery.  

If we choose the effective resistance edge-weights in \Cref{alg-main}, then \Cref{thm-1} and \Cref{prop-1} imply that the error rate we have obtained is nearly minimax optimal.  

Furthermore, recall from \Cref{prop-2} that no algorithm is guaranteed to be consistent in the low signal-to-noise regime  $\kappa^2 \Delta \sigma^{-2}  \lesssim |\partial_r (\mathcal{S}^*)|$. On the other hand, letting the edge weights be the effective resistance weights, \Cref{thm-1} shows that consistent estimation is possible when  $\kappa^2 \Delta \sigma^{-2}  \gtrsim |\partial_r (\mathcal{S}^*)|$. Thus,  \Cref{thm-1} additionally reveals the existence of a phase transition in the space of model parameters, as in the high signal-to-noise ratio regime, consistent estimation is not only feasible but it can be done at a nearly minimax optimal rate.

\begin{remark}[Edge weights]
The edge-weights play an important role in all the results we have shown so far.  The two most commonly-used choices are: (1) the 0/1 weighting, i.e.~assigning unit weight to each edge, and (2) the effective resistance weighting.  Since $r(i, j) \leq 1 = \mathbf{1}\{(i, j) \in E\}$, for all $(i, j) \in E$, for the partition recovery task, our theory suggests that the  effective resistance weighting should be preferred to 0/1 weighting.  Having said this, when the sizes of $V$ and $E$ are moderately large, it is much more practical to directly adopt the 0/1 weighting rather than calculating the effective resistance.  This is also what we do in \Cref{sec-numeric}.
\end{remark}

\subsection{General \texorpdfstring{$K^*$}{TEXT}}	

In \Cref{sec-K=2} we only consider the case when $K^* = 2$ and is known.  There are ample applications where such cases are interesting and we refer to the Introduction of \cite{arias2011detection}.  Despite the popularity of the case $K^* = 2$, there are of course many interesting situations where $K^* \neq 2$.  In  chain graphs, an $\ell_0$-penalised method is shown to be optimal in change point localisation for general $K^*$, where $K^*$ is seen as the number of change points plus one, and is even allowed to diverge with the number of nodes $n$ \citep[e.g.][]{wang2020univariate, verzelen2020optimal}.  This result has been further extended to $d$-dimensional lattice graphs in  \cite{padilla2021lattice}, where under stronger conditions, it is shown that a constrained $\ell_0$-penalised estimator (dyadic classification and regression trees, DCART) is able to handle general $K^*$, despite possessing a gap regarding  $K^*$ in relation to a minimax lower bound.  For general graphs, when the goal is de-noising, \cite{fan2018approximate} showed that an $\ell_0$-penalised estimator is optimal for general $K^*$.  For us, with general graphs and partition recovery purpose, we show that the phenomenon is very different. 

\subsubsection{The Case of \texorpdfstring{$K^* = 1$}{TEXT}}
The first result we show is when $K^* = 1$, with large probability, the output of \Cref{alg-main} is constant over the whole graph.  This claim is formalised in \Cref{prop-constant}, which is very similar to Theorem~3.5 in \cite{fan2018approximate}.

\begin{proposition}\label{prop-constant}
Let $\delta \leq \sigma/\sqrt{n}$, $\tau \leq \sigma^2$ and $\lambda = C_{\lambda} \sigma^2 \log\{w(|E|)\}$, where $C_{\lambda} > 0$ is an absolute constant.  Assume that $K^* = 1$.  Under \Cref{assume-model}, letting $\widehat{\mu}$ be an output of \Cref{alg-main} with $w \geq r$, it holds with probability at least $1 - \{w(|E|)\}^{-c}$ that $|\widehat{\mathcal{S}}| = 1$, where $\widehat{\mathcal{S}}$ is the partition induced by $\widehat{\mu}$.
\end{proposition}

\subsubsection{The Case of Unknown \texorpdfstring{$K^* > 1$}{TEXT}}
One may wish to say that an immediate consequence of \Cref{thm-1} and \Cref{prop-constant} is that when $K^* = 2$ but $K^*$ is unknown, \Cref{thm-1} still holds by first conducting \Cref{alg-main} on the whole graph, then repeatedly and separately on the resulting two pieces. (For completeness, we formalise this procedure in \Cref{alg-general}, with the notation that for any two partitions $\mathcal{A}$ and $\mathcal{B}$, $\mathcal{A} \cap \mathcal{B}$ is their refinement.)

\begin{algorithm}[ht]
	\begin{algorithmic}
	    \INPUT $G = (V, E)$, $|V| = n$, $\{Y_i, \, i \in V\}$, $\delta, \tau, \lambda, w$.
	    \State $\widetilde{\mu} \leftarrow 0^{\otimes n}$, $\mathrm{FLAG} \leftarrow 0$
	    \While{$\mathrm{FLAG} = 0$}
	        \State $\widehat{\mathcal{S}}, \widetilde{\mathcal{S}} \leftarrow $ the induced partition of $\widetilde{\mu}$ 
	        \For{each $S \in \widehat{\mathcal{S}}$}
	            \State $\mu_1 \leftarrow \mathrm{Potts}(G_S, Y_{G_S}, \delta, \tau, \lambda, w)$  \hspace{-1cm}\Comment{\Cref{alg-main}}
	            \State $\mathcal{S}_1 \leftarrow $ the induced partition of $\mu_1$
	            \State $\widehat{\mathcal{S}} \leftarrow \widehat{\mathcal{S}} \cap \mathcal{S}_1$
	            \State $\widetilde{\mu}_S \leftarrow \mu_1$ \Comment{$\widetilde{\mu}_S$ is the sub-vector of $\widetilde{\mu}$ on $S$}
	        \EndFor
	        \State $\mathrm{FLAG} = \mathbf{1}\{\widehat{\mathcal{S}} = \widetilde{\mathcal{S}}\}$
	    \EndWhile
	    \OUTPUT $\widehat{\mathcal{S}}$
		\caption{Variant of \Cref{alg-main}} \label{alg-general}
	\end{algorithmic}
\end{algorithm} 

Unfortunately, such result can only hold under a  much stronger assumption and with a much worse rate.  A direct consequence of \Cref{thm-1} is the following.  

Assume that there exists a constant $C_{\mathrm{SNR}} > 0$ such that 
    \begin{equation}\label{eq-snr-strong}
        \kappa^2 \Delta \geq C_{\mathrm{SNR}} \sigma^2 n \log(n) \zeta_n,
    \end{equation}
    where $\zeta_n$ is any arbitrarily diverging sequence, as $n$ grows unbounded.  Let $\delta \leq \sigma/\sqrt{n}$, $\tau \leq \sigma^2$ and $\lambda = C_{\lambda} \sigma^2 \log\{w(|E|)\}$, where $C_{\lambda} > 0$ is an absolute constant.  Assume that $K^* = 2$ but $K^*$ is unknown, and assume that \Cref{assume-model} holds.  Letting $\widehat{\mathcal{S}}$ be an output of \Cref{alg-general}, with $w \geq r$, it holds with probability at least $1 - \{w(|E|)\}^{-c}$ that 
    \begin{equation}\label{eq-strong-results}
        d_{\mathrm{H}}(\widehat{\mathcal{S}}, \mathcal{S}^*) \leq C \kappa^2 \sigma^{-2} n \log(n) \zeta_n.
    \end{equation}

Comparing \Cref{thm-1} and the aforementioned, we see that when $K^*$ is unknown, the corresponding rates in both the signal-to-noise ratio condition and the estimation error bound, jump from $|\partial_w(\mathcal{S})|$ to $n$.  Essentially, this is due to the complexity of the graphs and we elaborate as follows.

The proof of \Cref{thm-1} is built upon the large probability event defined in \eqref{eq-large-prob-event-fan}, and \eqref{eq-strong-results} is due to the event
\[
\left\{\left||A|^{-1/2} \sum_{i \in A} \varepsilon_i \right| \leq C n \sigma, \quad \forall A \in 2^V\right\},
\]
where $2^V$ is the power set of $V$. 
On both events, the difference between the sample and population quantities are upper bounded.  The large probability events are constructed based on sub-Gaussian concentration inequalities and a union bound argument.  In \Cref{thm-1}, since $K^* = 2$ and $K^*$ is assumed to be known, and due to the assumption that $G$ is a connected graph, the union bound argument is based on a certain spanning tree.  This is shown in Lemma~B.2 in \cite{fan2018approximate}.  In \Cref{alg-general}, due to the different layers of splits, at each layer, one works based on a random sub-graph yielded by the previous layer optimisation.  To take this randomness into consideration, one therefore has to consider the complexity of the whole graph.  

To further understand how the complexity kicks in, we note that in the change point localisation literature, i.e.~when $G$ is a chain graph, there is no such dramatic change in rates for general $K^*$ cases \citep[e.g.][]{wang2020univariate}.  This is due to the fact that there is one and only one spanning tree in a chain graph, and any spanning tree of any sub-graph is a sub-tree of a spanning tree of the whole graph.  For a general graph, it is possible to partition a sub-graph without cutting through a spanning tree of the whole graph.  Therefore, when applying a union bound argument over random sub-graphs, one has to consider all possible spanning trees.  The number of all spanning trees is  exponential in $n$, which holds even when the graph is a $k$-regular graph, with $k$ a fixed constant \citep[e.g.][]{alon1990number}, or when the graph is a lattice \citep[e.g.][]{shrock2000spanning}.

Based on the above discussions, one may be led to  conjecture that for general $K^*$, \eqref{eq-snr-strong} and \eqref{eq-strong-results} are in fact optimal.  If true, this finding would be rather interesting -- if the signal-to-noise ratio has to be as large as required in \eqref{eq-snr-strong}, then for any edge $(i, j) \in E$, one can simply cut this edge if and only if $|Y_i - Y_j| \gtrsim \sigma \log(n)$ and optimally recover the partition.
  
\subsection[]{Comparisons with \cite{fan2018approximate}}\label{sec-fan}

\cite{fan2018approximate} is the closest-related literature and has heavily inspired our work.  We therefore provide a thorough comparison with \cite{fan2018approximate}.  Though the set-up is identical, the contributions of the present paper are markedly different in both the goals and technical features than those in \cite{fan2018approximate}. Below, we highlight how our work relates to \cite{fan2018approximate}. 
(a) Task. \cite{fan2018approximate} investigate de-noising or signal estimation of a piece-wise constant signal over a graph while  the present paper aims to estimate the boundary of the partition supporting the constancy regions of the signal - a problem that, somewhat surprisingly, has never been thoroughly studied in the literature.  Despite their apparent similarities, de-noising and partition recovery  are fundamentally different tasks, with different metrics and challenges. For the specific case of the chain graph, the two tasks have led to distinct lines of research and very different results.  In particular, the partition recovery problem is effectively equivalent to change point localisation in a time series.
(b) Optimality and assumptions.  Both \cite{fan2018approximate} and ours provided minimax optimal results, supported by minimax lower bounds.  The lower bound proof in \cite{fan2018approximate} is an application of Theorem 1(b) in \cite{raskutti2011minimax}, which is based on the metric entropy of $\ell_q$ balls.  We have three lower bound results, which rely on graph theory. 
(c) Algorithms.  \Cref{alg-main} and that in \cite{fan2018approximate} share the same core procedure, which is the $\alpha$-expansion algorithm \cite{boykov2001fast}.  Our \Cref{alg-main} is indeed a straightforward adaptation of that in \cite{fan2018approximate} by simply imposing that only one split is performed.  
(d) The proofs of the theoretical guarantees of \Cref{alg-main}.  The aforementioned differences between the problems imply that the structure and key technical aspects of the proofs are completely different.  Both sets of proofs rely on a large probability event put forward in \cite{fan2018approximate}, which we directly cite. 
(e) Other aspects.  Due to the different nature of the problem, \cite{fan2018approximate} is able to handle a growing number of constancy pieces, while such case becomes challenging and yet to be completely understood in the partition recovery case.  However, through comparisons with change point detection in chain graphs and abnormal cluster detection in general graphs, we provide more insights on the role of graph sparsity in such problems.

\section{NUMERICAL EXPERIMENTS}\label{sec-numeric}

We now proceed to evaluate the empirical performance of the proposed method, with four cases described below.  In each case we consider data in the 2D grid graph $G$ with $\sqrt{n} \times \sqrt{n}$ nodes,  where $K^* = 2$ is known, $n  \in \{64^2,128^2\}$,  $\sigma = 1$ and $\kappa \in \{1,2\}$.  We consider three competing methods:  \Cref{alg-main}, a variant of the dyadic CART \citep[DCART,][]{donoho1997cart} and the edge lasso \citep[EL,][]{sharpnack2012sparsistency}. For each method and case we report the median Hausdorff distance based on 50 repetitions.

For the implementation of \Cref{alg-main}, we set $\delta = 1/60$ and choose $\lambda$  from the set of candidates $\{10^{-2},10^{-1},10^0,10^1,10^2\}$  via the Bayesian information criterion (BIC). Specifically, we first estimate $\sigma^2$ with $\hat{\sigma}^2$ given as $\hat{\sigma}^2 = (n -1)^{-1}\sum_{i=1}^{n-1} (Y_{ p(i) }- Y_{p(i+1)})^2$,
where $p(1), \ldots, p(n)$ form a connected path of the grid graph.  With this in hand, for any $\lambda$ and its corresponding $\tilde{\mu}(\lambda)$ we calculate the score
\begin{equation} \label{eqn:bic}
	\mathrm{BIC}_{\lambda} = \sum_{i=1}^{n} \{Y_i - \tilde{\mu}_i(\lambda)\}^2 + \hat{\sigma}^2 v_{\lambda} \log(n),
\end{equation}
where  $v_{\lambda}$  is the number of connected components induced by $\tilde{\mu}(\lambda)$ in $G$. We then choose the value of  $\lambda$ that minimizes  $\mathrm{BIC}_{\lambda}$ and report the performance of  \Cref{alg-main} based on this choice.
The implementation of  \Cref{alg-main} is done in Matlab using the package   \texttt{gco-v3.0} \citep{boykov2001fast}.

\begin{table}[h]
\caption{Median Hausdorff distances of 50 repetitions under different cases.  DCART: dyadic CART; EL: edge Lasso; $\kappa$: jump size; and $n$: number of nodes.} \label{tab1}
\begin{center}
\begin{tabular}{llllll}
\textbf{CASE}  & $\kappa$ & $\sqrt{n}$ & \textbf{ALG.~1} & \textbf{DCART} & \textbf{EL} \\
\hline \\
$1$            &       $1.0$ &    128       &  \textbf{56} &  115 &922\\  
			$1$            &       $2.0$ &      128    &   33    &86& \textbf{21}\\  
			$1$            &       $1.0$ &        64   &   \textbf{29}      & 152 & 1630 \\  
              $1$            &       $2.0$ &       64   &  12   &34& \textbf{10}\\  [3pt]
			$2$            &       $1.0$ &    128       &  \textbf{122}   & 230 & 516\\  
            $2$            &       $2.0$ &      128    &   42  &160&  \textbf{32}\\  
            $2$            &       $1.0$ &        64   & \textbf{89}   &   294 &1601 \\  
            $2$            &       $2.0$ &       64   &  24  &63&\textbf{15}\\  [3pt]
  			$3$            &       $1.0$ &    128       &   36  &  63&\textbf{7}\\  
             $3$            &       $2.0$ &      128    &  \textbf{1}&  15&\textbf{1}\\  
              $3$            &       $1.0$ &        64   &     \textbf{20}&   47 & 1392\\  
             $3$            &       $2.0$ &       64   &   10  &15&\textbf{1}\\   [3pt]
          $4$            &       $1.0$ &    128       & \textbf{514}   &   604  &6669\\  
          $4$            &       $2.0$ &      128    &    \textbf{83} &199&671\\  
          $4$            &       $1.0$ &        64   &     \textbf{387}&   424&  1836\\  
        $4$            &       $2.0$ &       64   &   \textbf{49}  &108&226\\ 
\end{tabular}
\end{center}
\end{table}

For DCART, we first compute their values with the penalty $\lambda \in \{10^{-1+ 4j/19 }:\, j  = 0, 1,\ldots, 19 \}$.  We then select $\lambda$ that minimises the BIC in \eqref{eqn:bic}, replacing  $\tilde{\mu}(\lambda)$ with the corresponding  DCART  estimator.  This produces an estimator $\hat{\mu}$. However, assuming that $K^* = 2$ is known, we let $\Lambda$ be the set of unique values of  $\hat{\mu}$ and perform $k$-means  clustering on the elements of $\Lambda$  setting the number of clusters  to  two. Let  $C_1$ and  $C_2$  be the centres obtained from the $k$-means clustering.  Define $\tilde{\mu}_i = C_1$, if $\hat{\mu}_i$ is assigned to the centre $C_1$; and $\tilde{\mu}_i = C_2$, otherwise.  The final estimator is the partition induced by $\tilde{\mu}$.  The computation for the resulting procedure was done in R.

For EL, we proceed as we do with DCART employing the $k$-means post-processing and BIC model selection. The choices of the penalty parameter are $\{10^{-2+6j/19}:\, j \in \{1,\ldots, 19\}  \}$.

The specific cases considered as as follows. In each case, we generate data as $Y_{i,j} = \mu^*_{i,j} + \varepsilon_{i,j}$, where  $\varepsilon_{i,j} \overset{\mathrm{i.i.d.}}{\sim} \mathcal{N}(0,1)$, $i,j = 1,\ldots, \sqrt{n}$ and $\mu^*$ is specified below.

\begin{figure}[t!]
	\begin{center}%
	\hspace{-0.05in}
	\includegraphics[width=0.9in,height=1in]{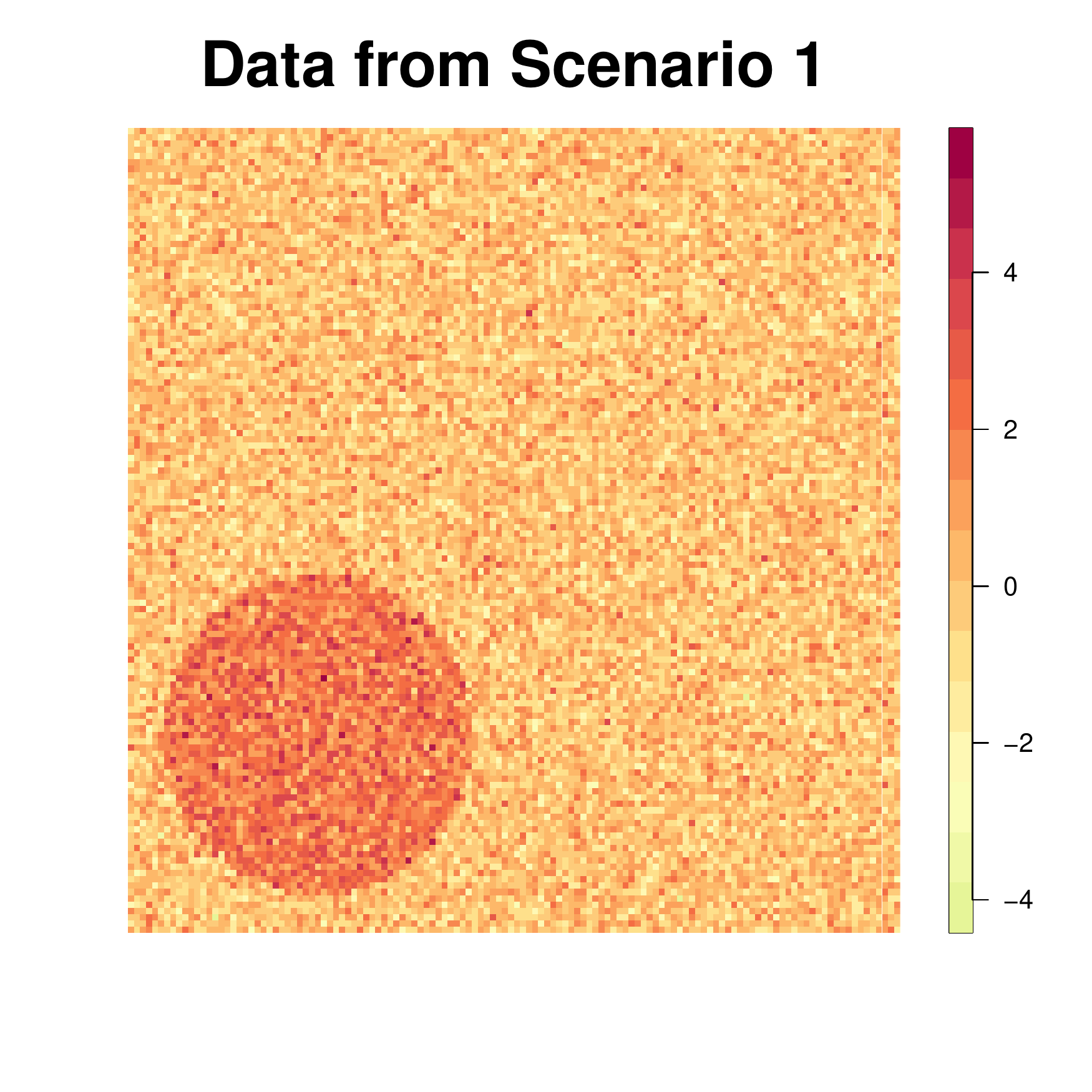} \hspace{-0.2in}
	\includegraphics[width=0.9in,height=1in]{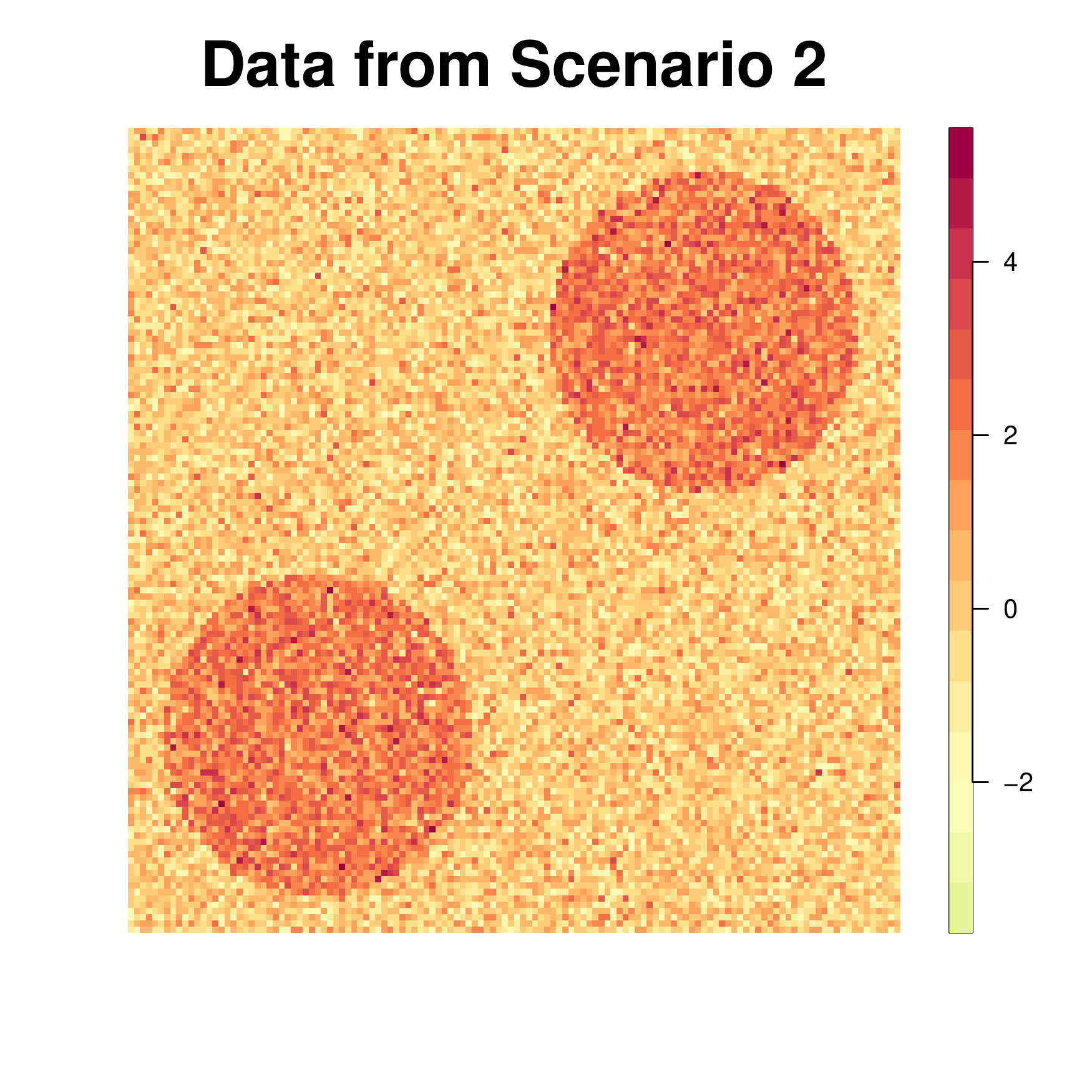}
	\hspace{-0.2in}
	\includegraphics[width=0.9in,height=1in]{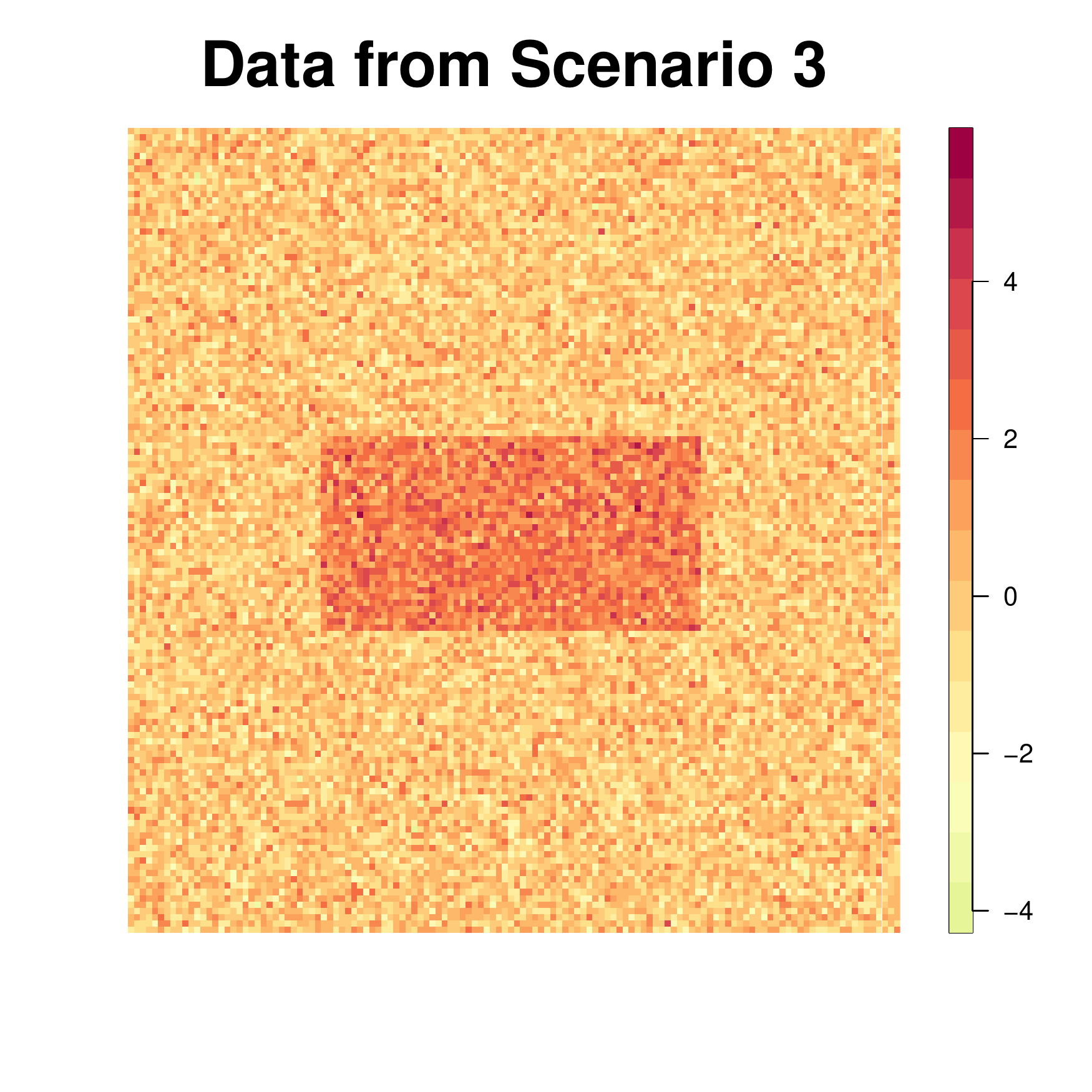}
	\hspace{-0.2in}
	\includegraphics[width=0.9in,height=1in]{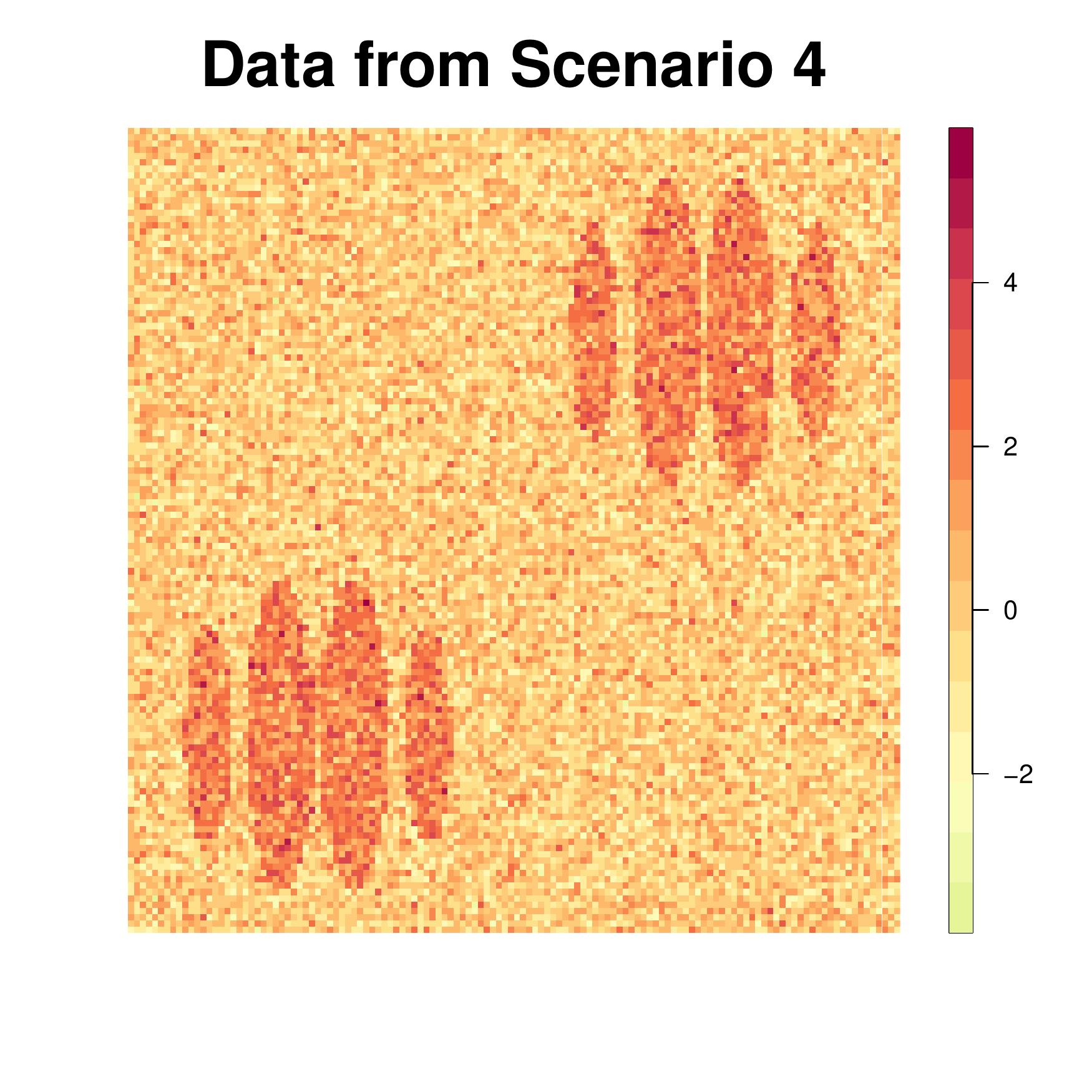}\\
		\hspace{-0.13in}
	\includegraphics[width=0.77in,height=1in]{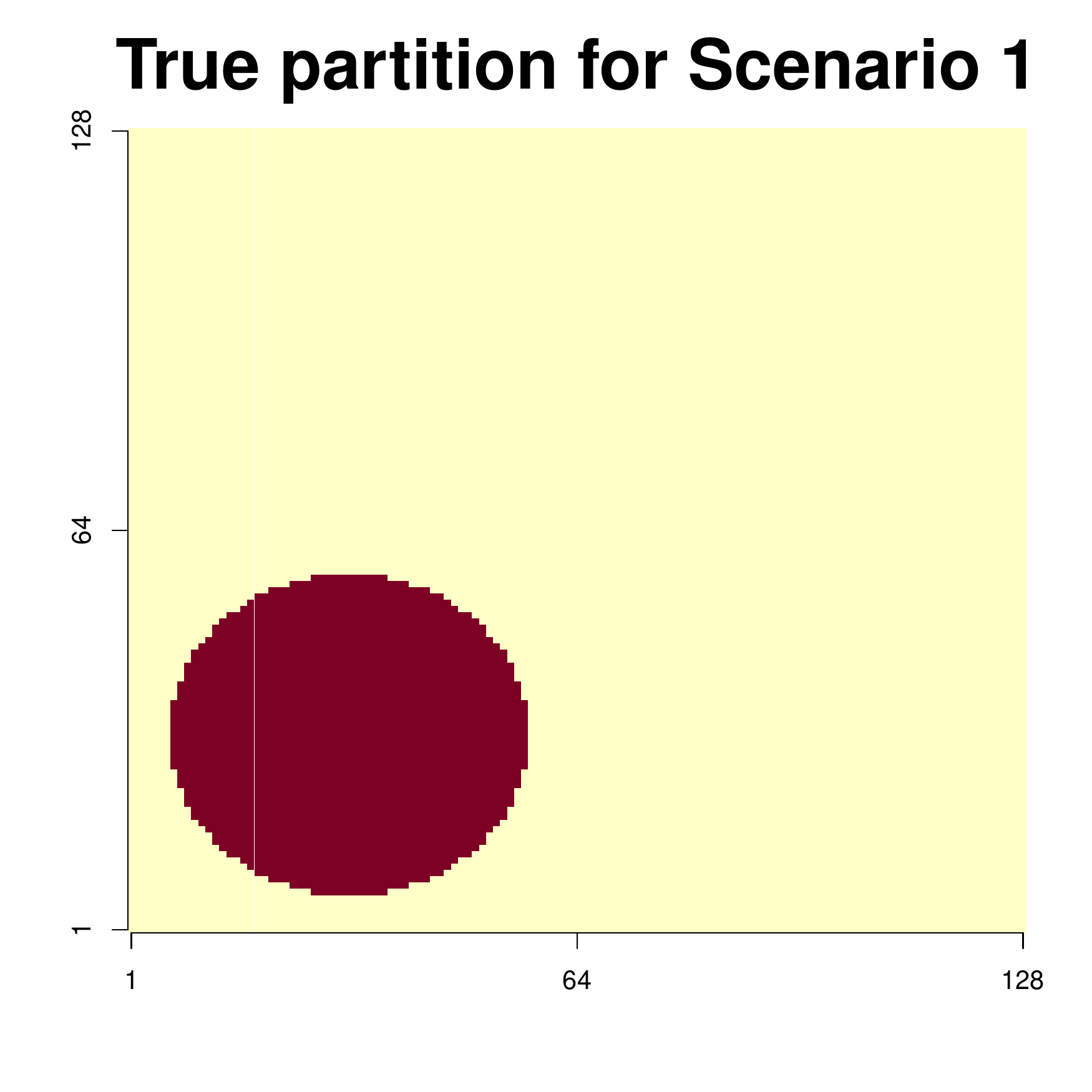}
	\hspace{-0.065in}
	\includegraphics[width=0.77in,height=1in]{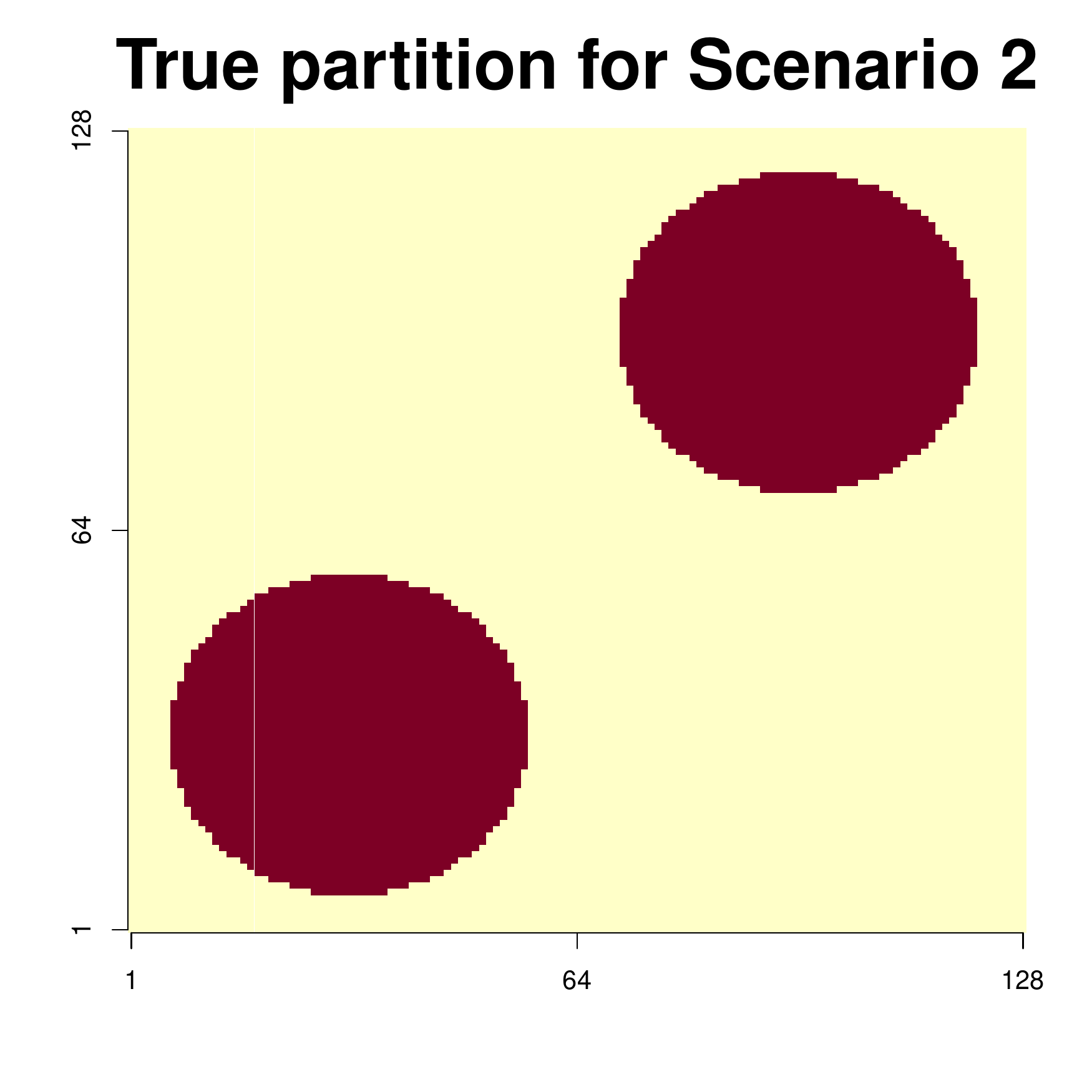}
		\hspace{-0.078in}
	\includegraphics[width=0.77in,height=1in]{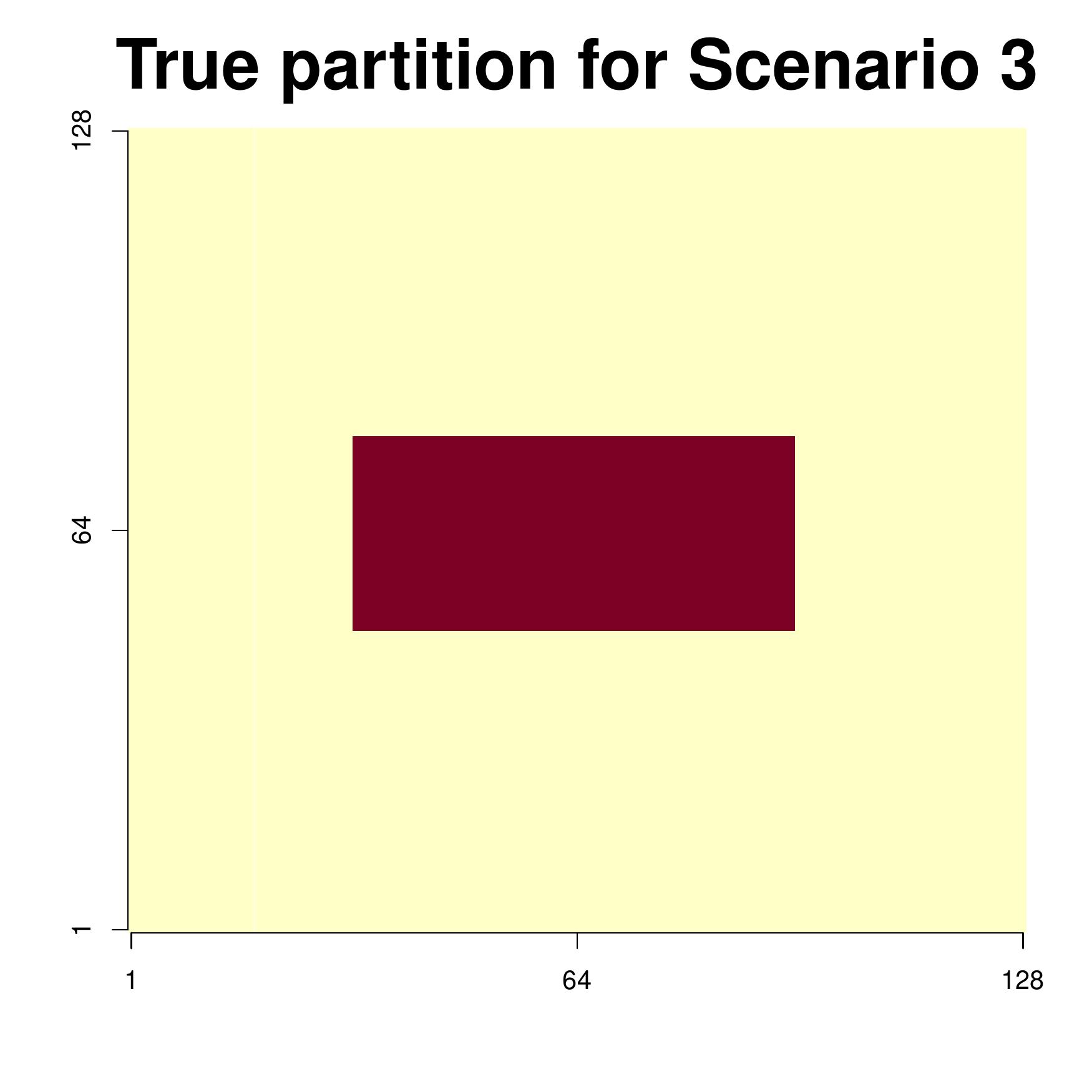}
		\hspace{-0.09in}
	\includegraphics[width=0.77in,height=1in]{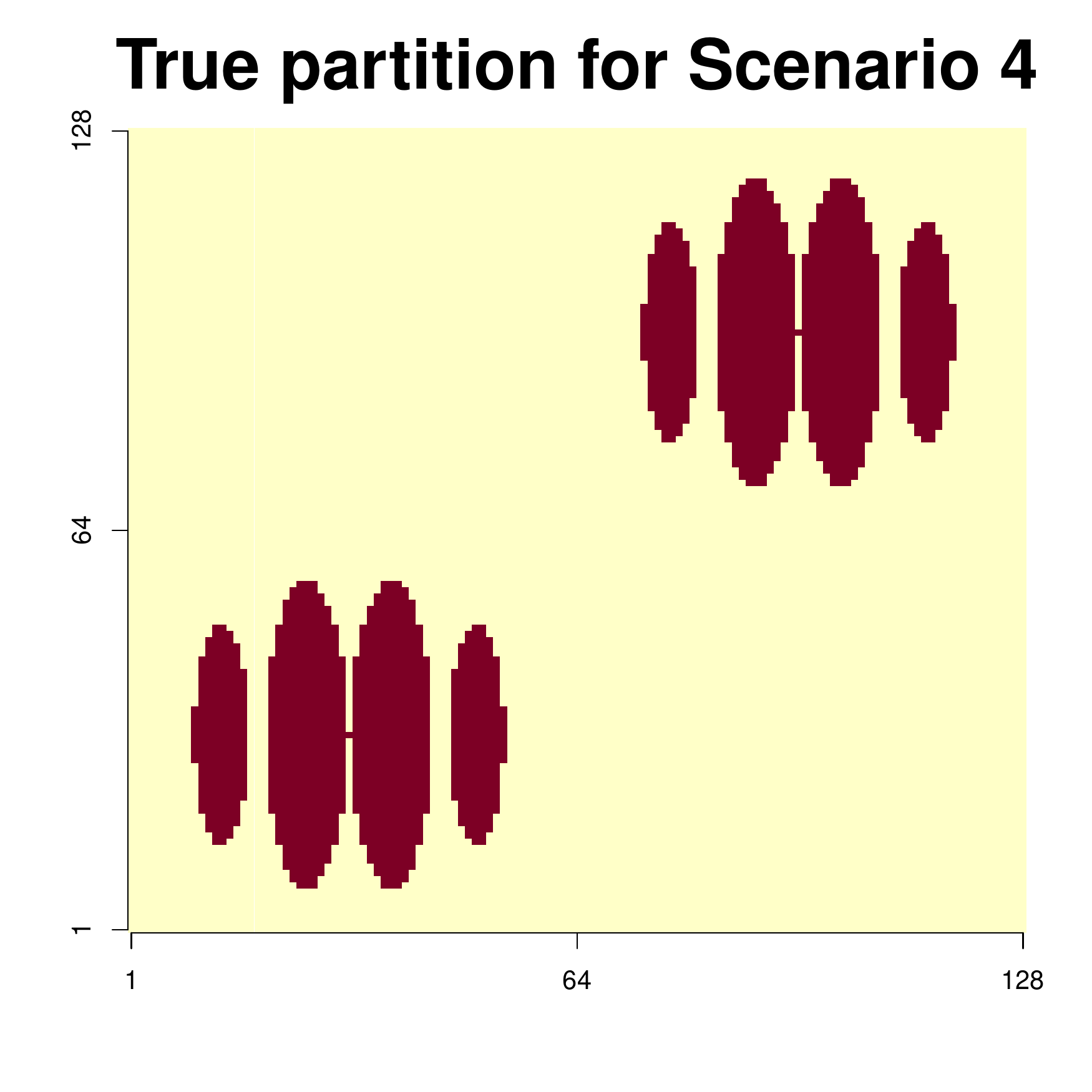} \\
		\hspace{-0.13in}
	\includegraphics[width=0.77in,height=1in]{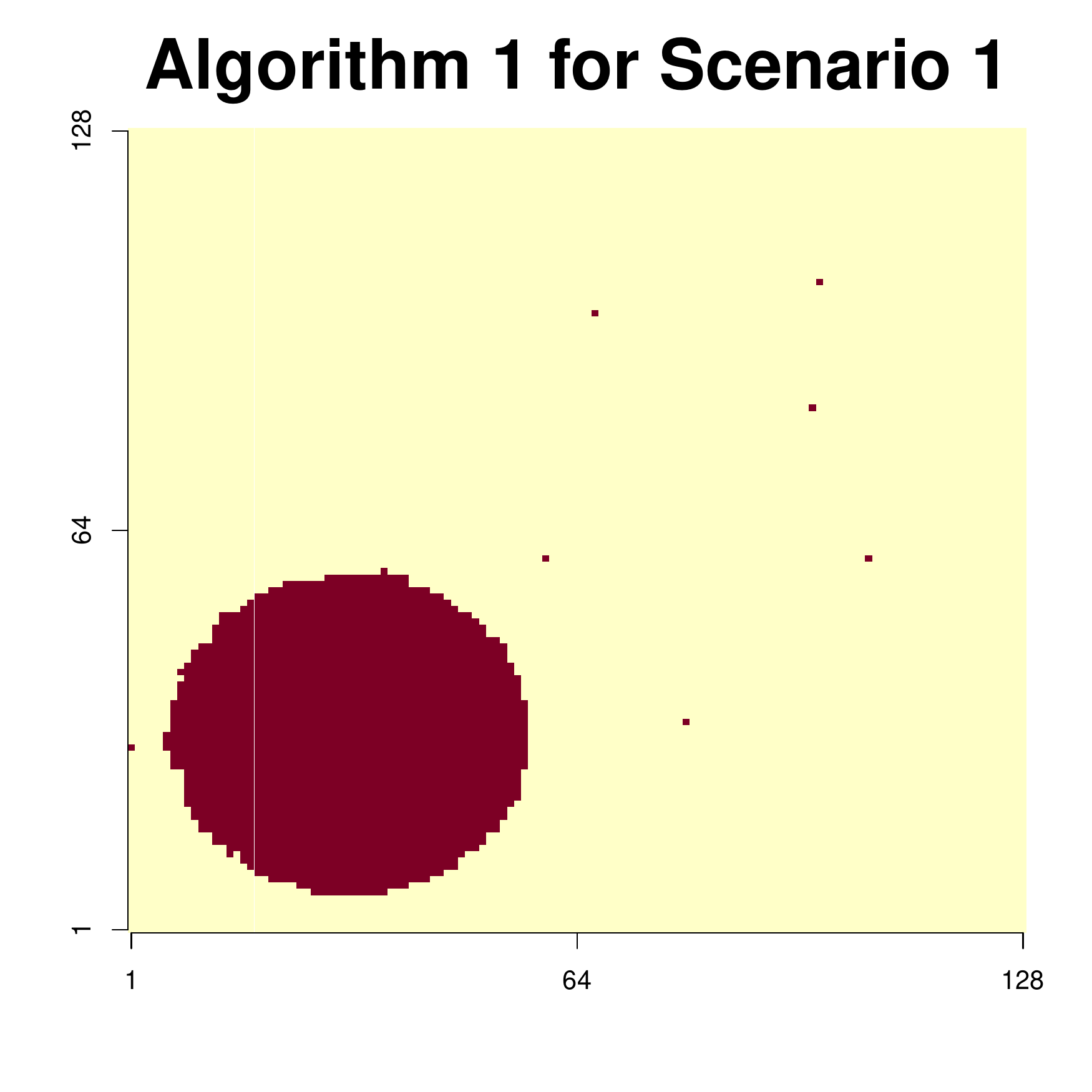} 
		\hspace{-0.065in}
	\includegraphics[width=0.77in,height=1in]{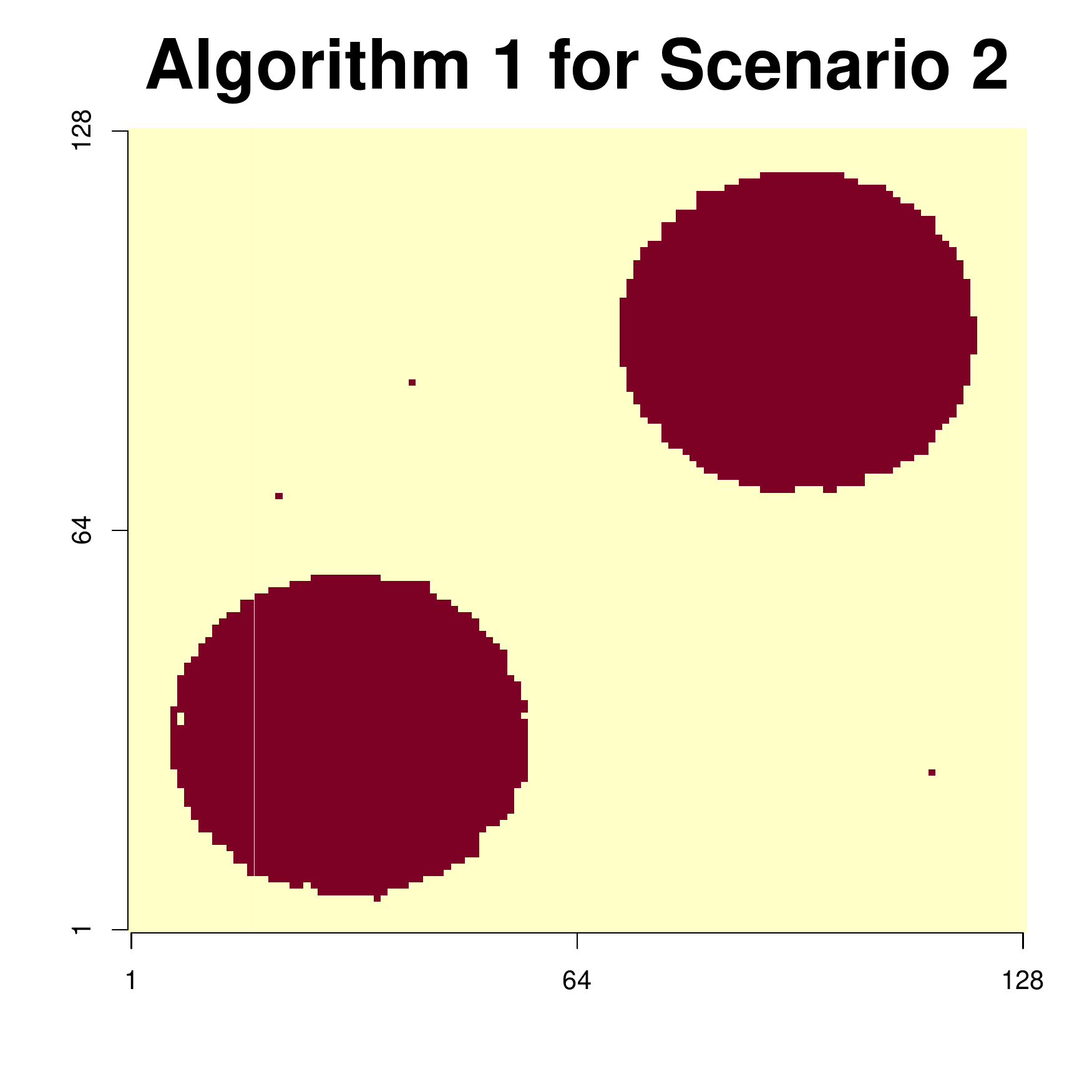}
		\hspace{-0.078in}
	\includegraphics[width=0.77in,height=1in]{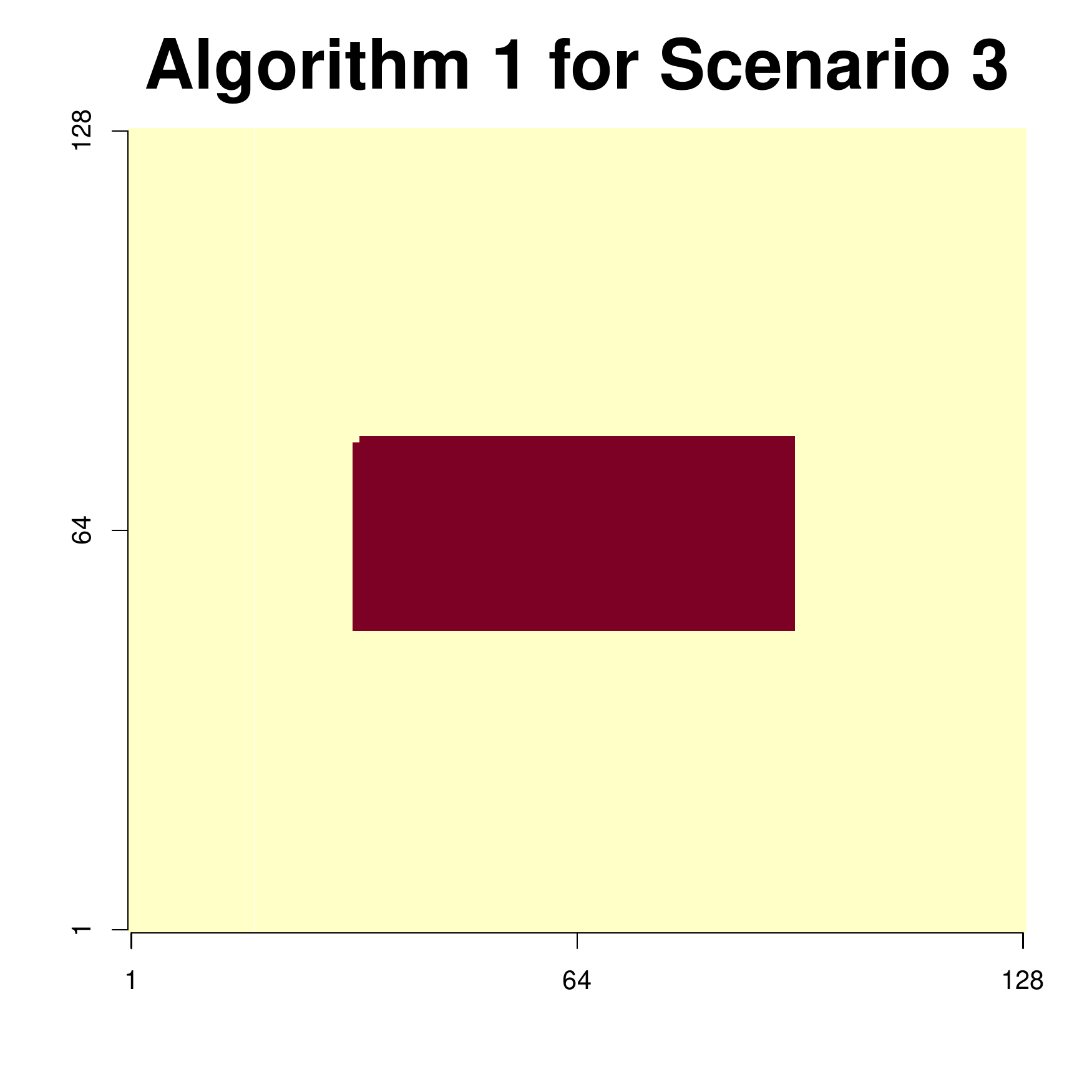} 
		\hspace{-0.09in}
	\includegraphics[width=0.77in,height=1in]{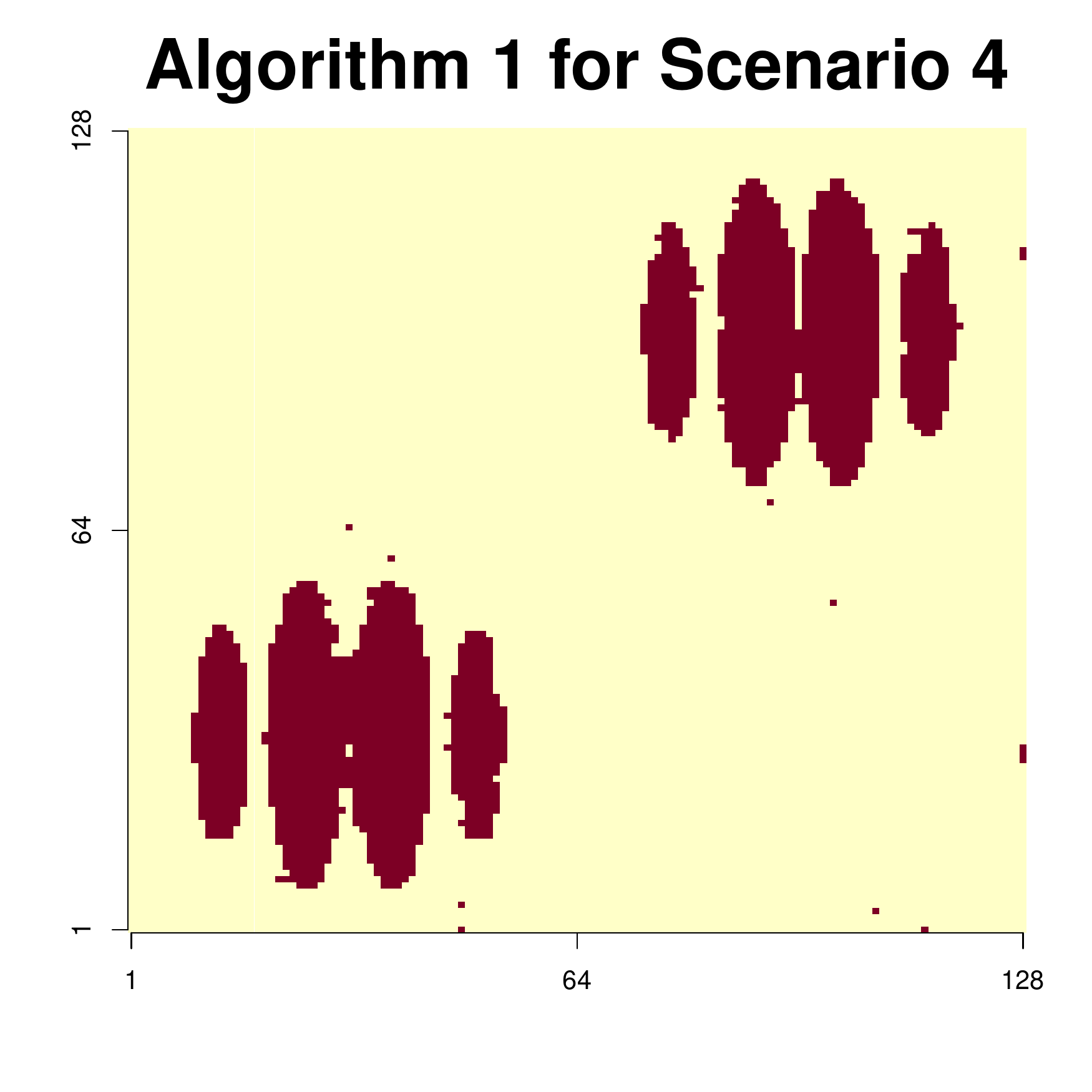} \\
		\hspace{-0.13in}
	\includegraphics[width=0.77in,height=1in]{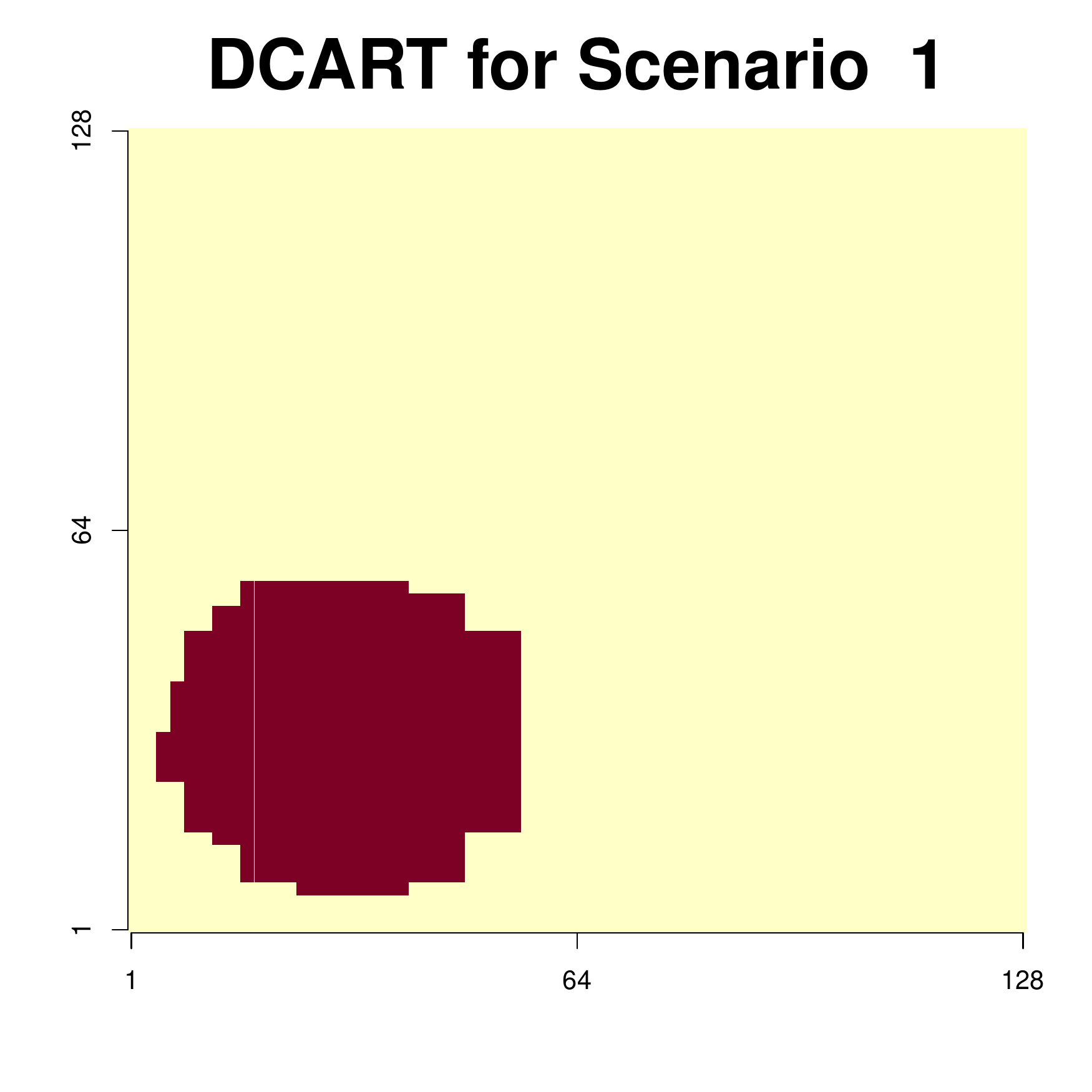} 
			\hspace{-0.065in}
	\includegraphics[width=0.77in,height=1in]{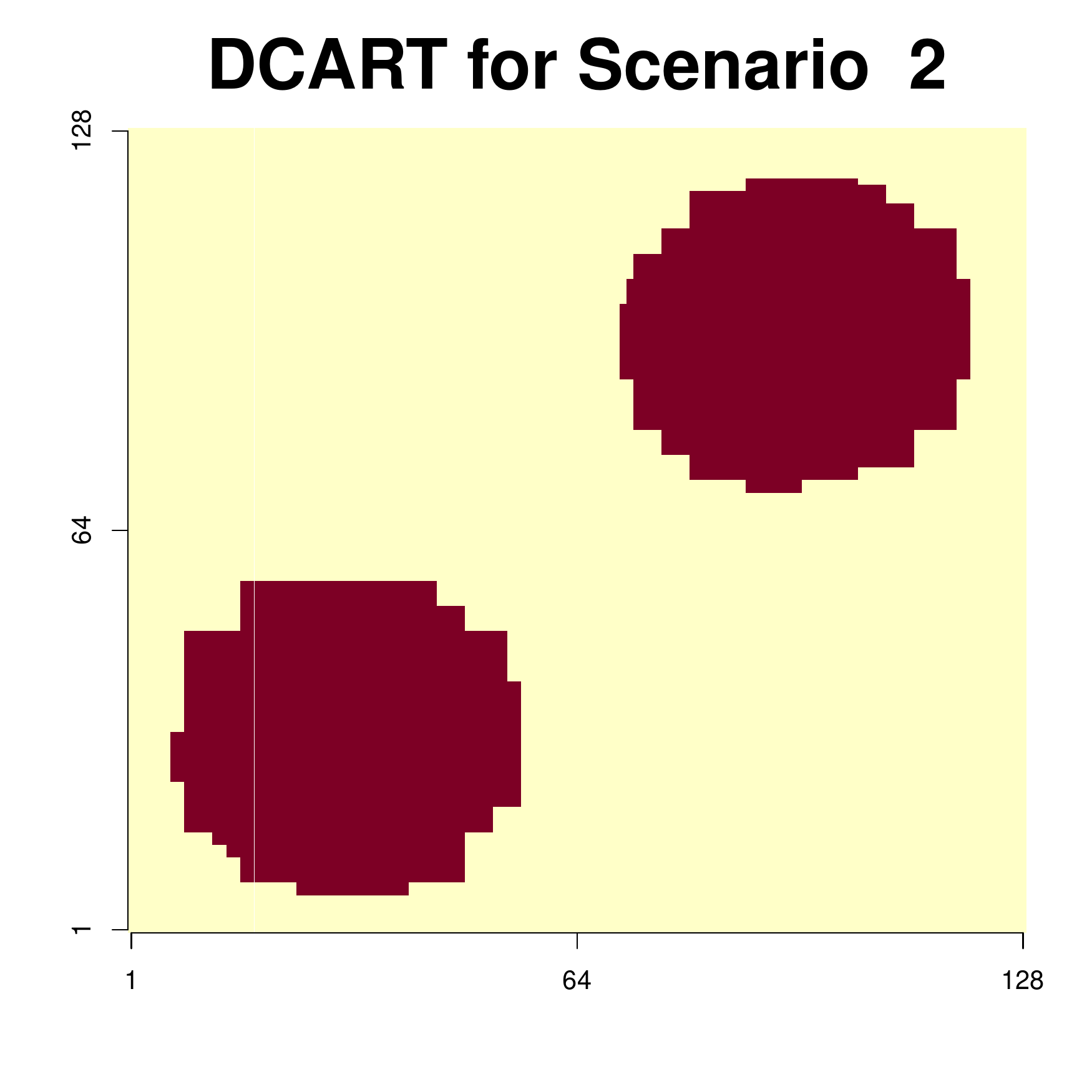}
			\hspace{-0.078in}
	\includegraphics[width=0.77in,height=1in]{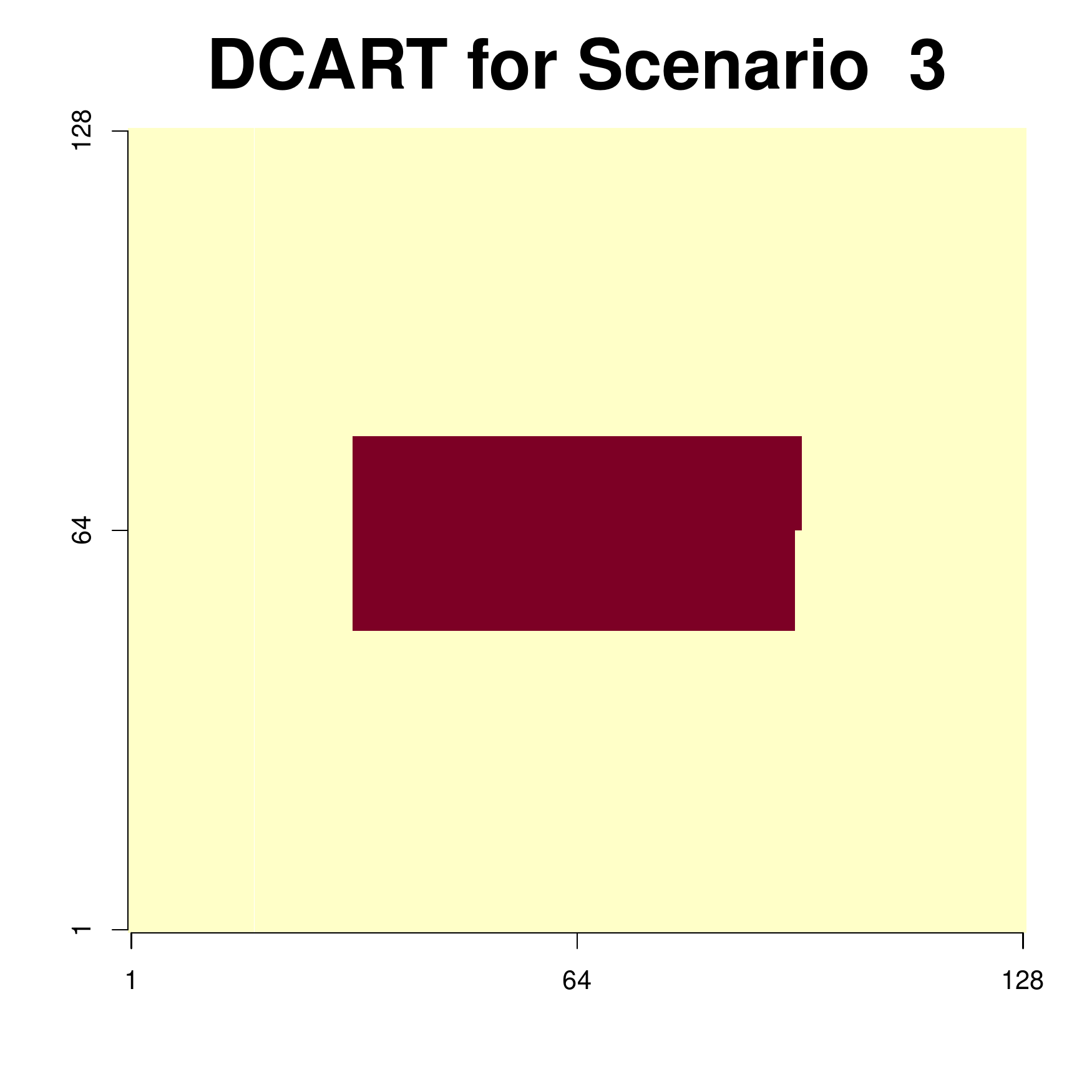}
		\hspace{-0.09in}
	\includegraphics[width=0.77in,height=1in]{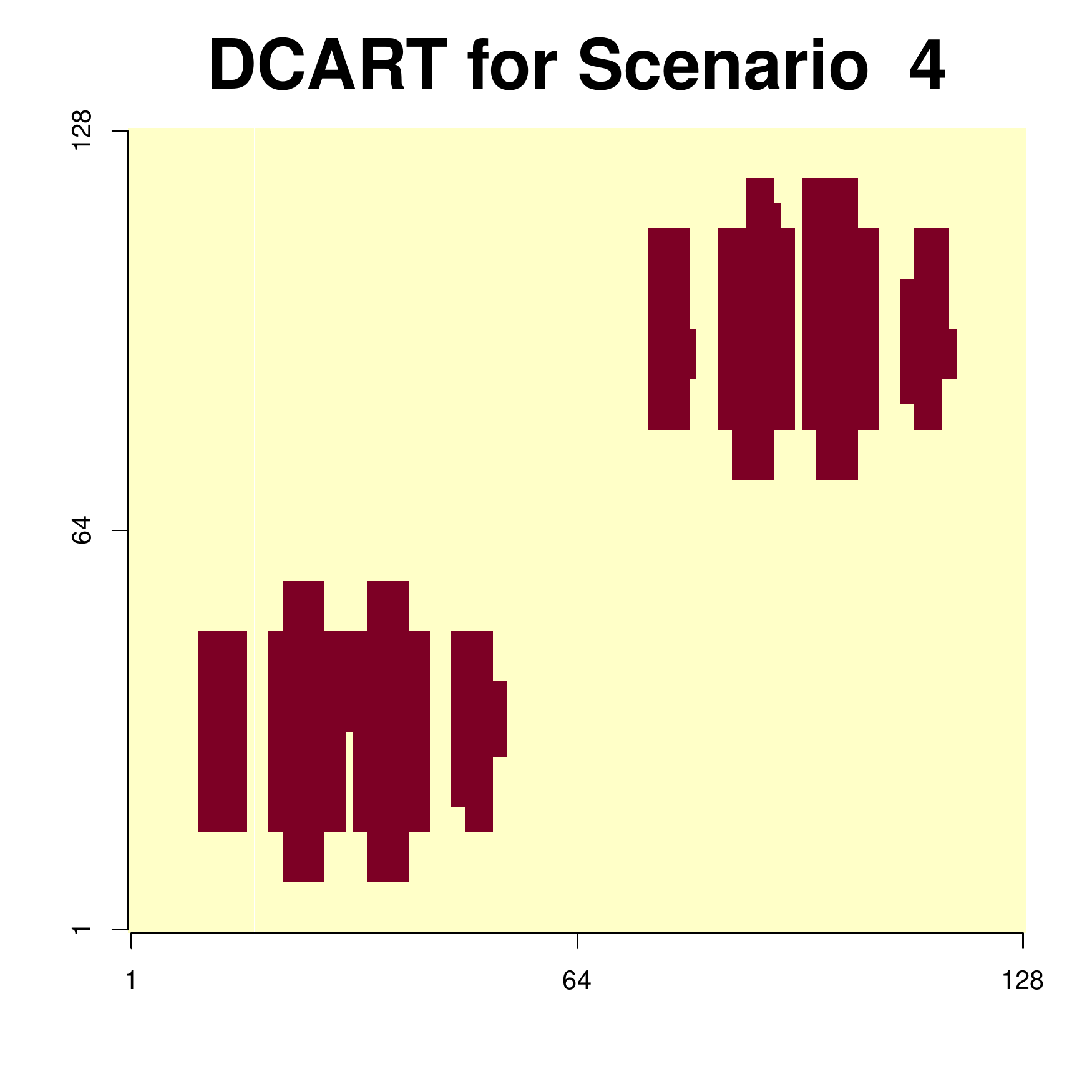}\\
		\hspace{-0.13in}
    \includegraphics[width=0.77in,height=1in]{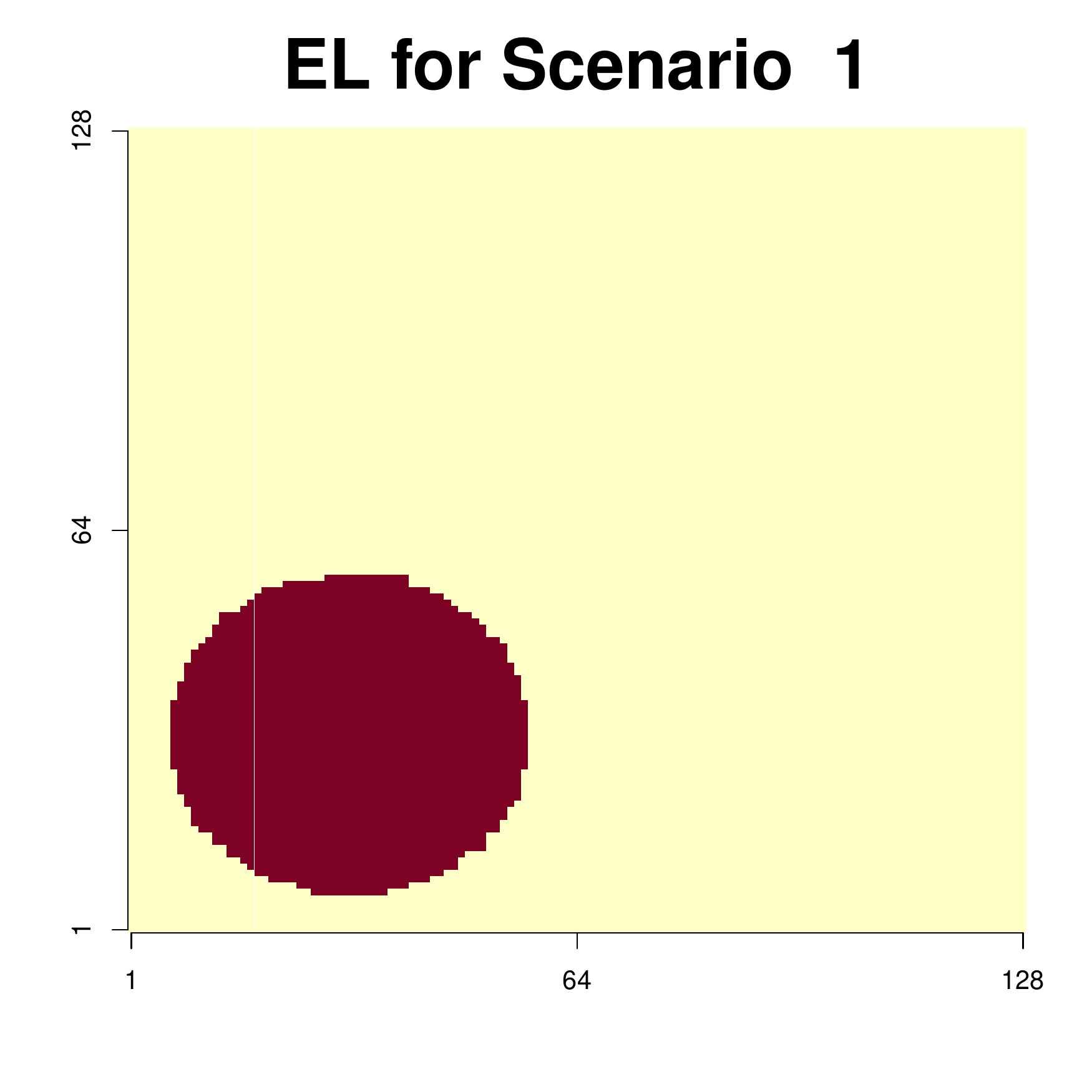} 
\hspace{-0.065in}    
    \includegraphics[width=0.77in,height=1in]{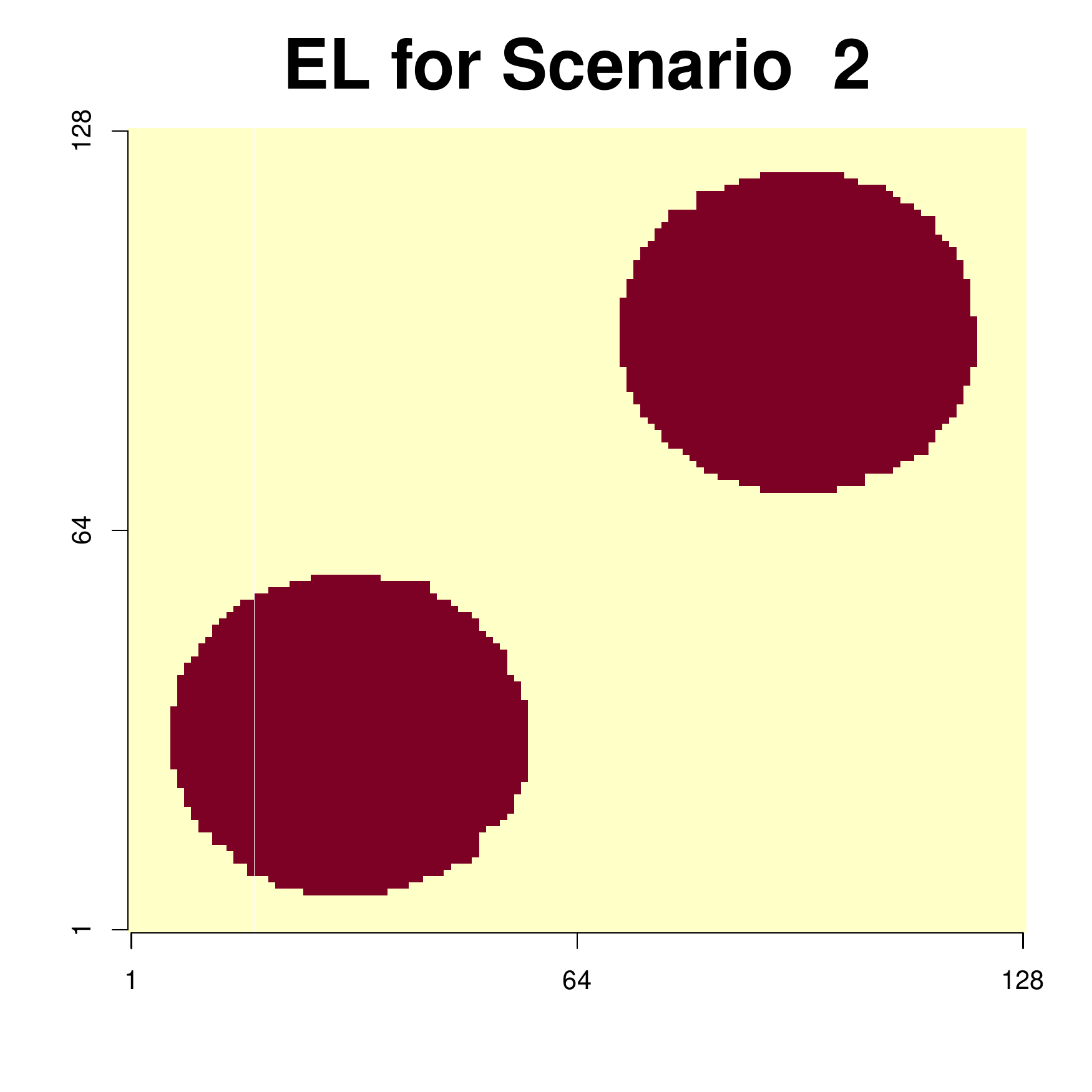} 
	\hspace{-0.078in}    
    \includegraphics[width=0.77in,height=1in]{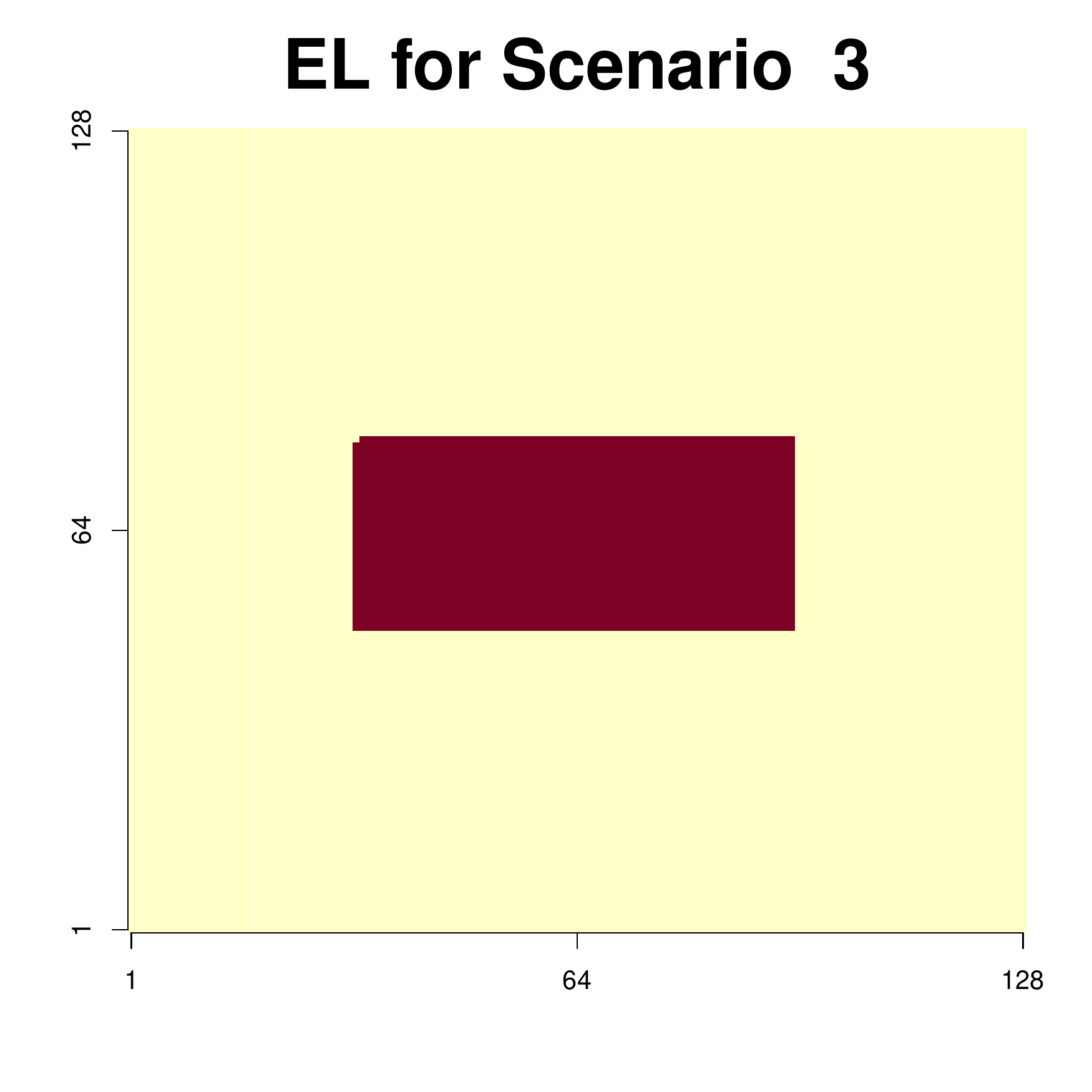} 
    		\hspace{-0.09in}
    \includegraphics[width=0.77in,height=1in]{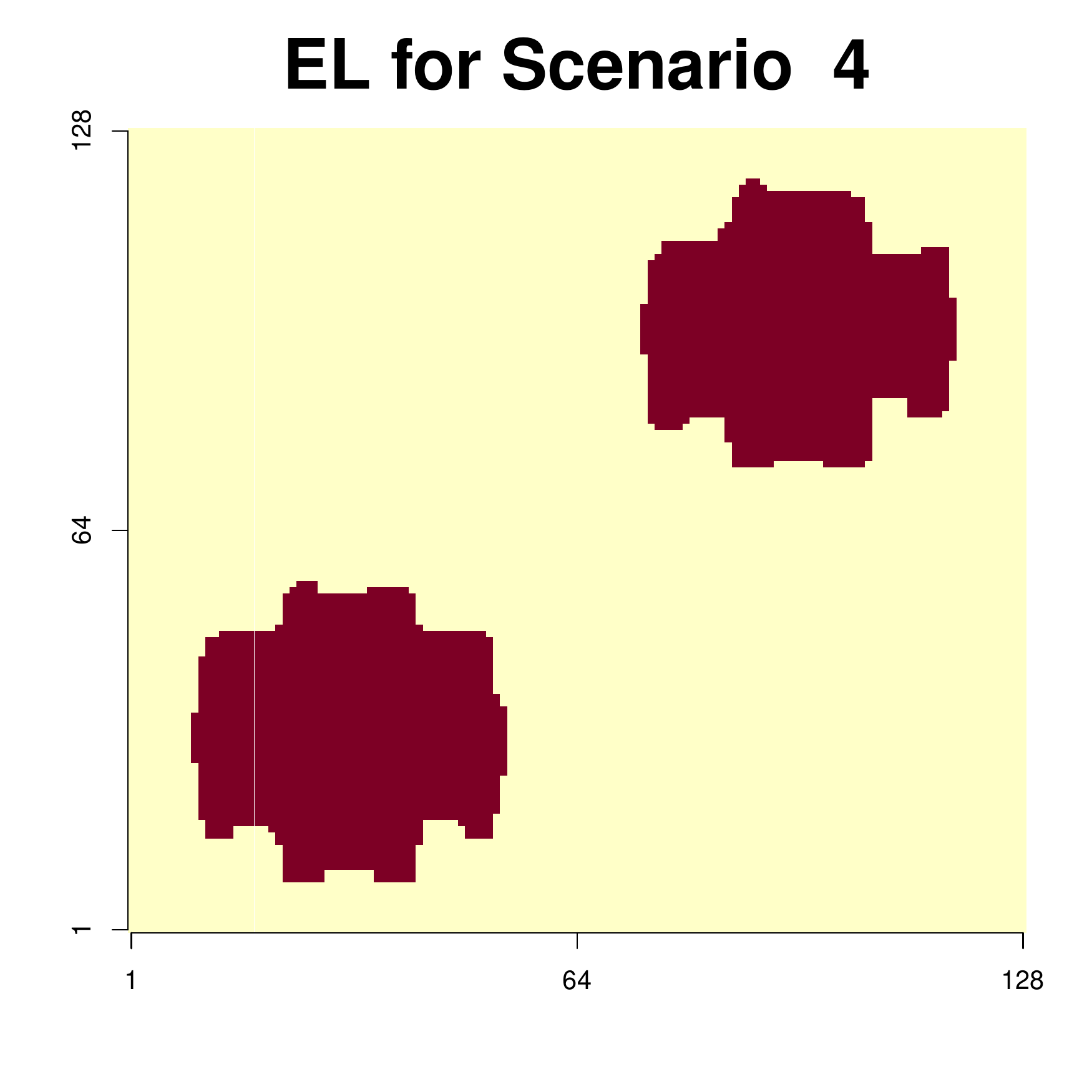} 
    \caption{\label{fig2} From top to bottom: instances of data, signal patterns, estimators of \Cref{alg-main}, dyadic CART (DCART) and edge lasso (EL).  From left to right: Cases 1, 2, 3 and 4, with $n = 128^2$ and $\kappa = 2$.}
	\end{center}
\end{figure}

\textbf{Case 1.} Let
\[
   \mu_{i,j}^* = \begin{cases}
   	\kappa,  & (i-  \sqrt{n}/4)^2 + (j-\sqrt{n}/4)^2< (\sqrt{n}/5)^2,\\
   	0, & \text{otherwise.}
   \end{cases}
\]

 \textbf{Case 2.} Let
 \[
 \mu_{i,j}^* = \begin{cases}
 	\kappa,  & (i-  \sqrt{n}/4)^2 + (j-\sqrt{n}/4)^2< (\sqrt{n}/5)^2\\
 	 & \text{or }(i-  3\sqrt{n}/4)^2 + (j-3\sqrt{n}/4)^2 \\
 	 & \hspace{0.5cm} < (\sqrt{n}/5)^2,\\
 	0, &\text{otherwise.}
 \end{cases}
 \]

 \textbf{Case 3.} Let
 \[
 \mu_{i,j}^*  = \begin{cases}
 	\kappa,  & \vert  i - \sqrt{n}/2\vert < \sqrt{n}/4 \\
 	 & \text{and }\vert  j - \sqrt{n}/2\vert < \sqrt{n}/4, \\
 	0, &\text{otherwise.}
 \end{cases}
 \]

 \textbf{Case 4.} Let
 \[
 \mu_{i,j}^*  = \begin{cases}
 	\kappa,  & (i-  \sqrt{n}/4)^2 + (j-\sqrt{n}/4)^2 \\
 	 &\,\,\,\,\,\,\, < \vert\mathrm{cos}(10\,\pi\,i/n)\vert (\sqrt{n}/5)^2\\
 	 & \text{or } (i-  3\sqrt{n}/4)^2 + (j-3\sqrt{n}/4)^2 \\
 	 &\,\,\,\,\,\,\, < \vert\mathrm{cos}(10\,\pi\,i/n)\vert (\sqrt{n}/5)^2,\\
 	0, &\text{otherwise.}
 \end{cases}
 \]

Visualisations of instances of generated data, true signals and different estimators are given in Figure \ref{fig2}.  We can see that the signal in Case 1 has two pieces, a circle and its complement.  Case 2 is based on two separated circles and their complement, but the signals on the two circles are the same.  The partition in Case 3 is a rectangle and its complement, making it the most attractive case for using the variant of DCART. The final case is perhaps the most challenging since the true signal has multiple pieces with boundaries that rapidly change in some regions.  Recall in \Cref{assume-model}, $S_k$'s are assumed to possess constant signal values, but are not necessarily connected.  Under this assumption, all four cases have $K^* = 2$.  Cases~2 and 4 show that although our method and theory are designed for $K^* = 2$, we enjoy great model complexity that we can handle cases with more than two connected constancy regions.

The results measured by the mean Hausdorff distances of 50 repetitions are shown in Table~\ref{tab1}, where we can see that \Cref{alg-main} performs better than the variant of DCART across all cases considered.  The EL can be competitive in cases where there is large signal-to-noise ratio but suffer greatly otherwise. 

\section{CONCLUSION}

In this paper, we study the partition recovery problem in general graphs, where nodes are associated with independent random variables.  We have shown that an $\ell_0$-penalised estimator is optimal when it is known that $K^* = 2$.  We have further derived a phase transition phenomenon, establishing the fundamental limits in the partition recovery problem in general graphs.  For general $K^*$, we have provided some seemingly unsatisfactory results, which we conjecture to be optimal and which prompt us discuss some crucial difference between the de-noising and partition recovery problems in general graphs. 

\subsubsection*{Acknowledgements}
All the authors thanks to the reviewers for constructive comments.  Yu is partially funded by EPSRC EP/V013432/1.  Madrid Padilla and Rinaldo are partially funded by NSF DMS 2015489.

\clearpage

\bibliography{ref}
\bibliographystyle{apalike}


\clearpage
\appendix

\thispagestyle{empty}

\onecolumn \makesupplementtitle

\section[]{PROOFS OF RESULTS IN SECTION~\ref{sec-lower-bound}}

Propositions~\ref{prop-2} and \ref{prop-1} can be directly shown by noticing that chain graphs are special cases of general graphs.  In chain graphs with $K^* = 2$, it holds that $|\partial_r(\mathcal{S}^*)| = 1$.  The results of Propositions~\ref{prop-2} and \ref{prop-1} follow directly from Lemmas~1 and 2 in \cite{wang2020univariate}.  In this section, we provide different proofs allowing $|\partial_r(\mathcal{S}^*)|$ to depend on the sample size $n$.

\begin{proof}[Proof of \Cref{prop-2}]
We prove the result by the Le Cam lemma \citep[e.g.~Lemma~1 in][]{yu1997assouad}.

\medskip
\noindent \textbf{Constructing distributions}
We consider a complete $n$-node graph $G = (V, E)$, where $V = \{1, \ldots, n\}$ and $E = \{(i, j), 1 \leq i < j \leq n\}$.

Given an absolute constant $c_1 > 0$, without loss of generality, we assume that $n/4$, $M = n/\{4c_1\log(n)\}$ and $c_1\log(n)$ are all positive integers.  For $l \in \{1, \ldots, M\}$, let $\widetilde{u}_l \in \mathbb{R}^n$ be such that the $i$the coordinate of $\widetilde{u}_l$, $i = 1, \ldots, n$, satisfies that
    \[
        \widetilde{u}_l(i) = \begin{cases}
            \sigma, & i \in \{(l-1)c_1 \log(n) + 1, \ldots, l c_1 \log(n)\}, \\
            0, & \mbox{otherwise}.
        \end{cases}
    \]
    Let $\widetilde{v}_l \in \mathbb{R}^n$ be such that $\widetilde{v}_l(i) = \widetilde{u}_l(n-i+1)$, $i = 1, \ldots, n$.  Let $\widetilde{P}_l$ and $\widetilde{Q}_l$ be multivariate Gaussian distributions $\mathcal{N}(\widetilde{u}_l, \sigma^2 I_n)$ and $\mathcal{N}(\widetilde{v}_l, \sigma^2 I_n)$, respectively and set
    \[
        \widetilde{P} = \frac{1}{M} \sum_{l = 1}^M \widetilde{P}_l \quad \mbox{and} \quad \widetilde{Q} = \frac{1}{M} \sum_{l = 1}^M \widetilde{Q}_l.
    \]

Note that for each $l \in \{1, \ldots, M\}$, $\widetilde{P}_l$ has two pieces, with $\Delta = c_1\log(n)$ and the smaller piece is contained in $\{1, \ldots, n/4\}$.  As a result, it holds that
    \begin{equation}\label{eq-proof-new-prop-2-1}
        \kappa \sqrt{\Delta} /\sigma = c_1 \log(n).
    \end{equation}
    In addition, since we are considering complete graphs, each edge $(i, j)$ has effective resistance weight
    \[
        r(i, j) = \frac{(n-1)n^{n-2}}{n(n-1)/2} = 2n^{n-3},
    \]
    following from the proof of \Cref{lem-com-resis-weight}.  Then
    \[
        |\partial_r(\mathcal{S})| = \frac{2n^{n-3}}{n^{n-2}} c_1 \log(n) \left\{n - c_1 \log(n)\right\} = \frac{2}{n}c_1 \log(n) \left\{n - c_1 \log(n)\right\}.
    \]
    Provided that $n \geq 2c_1 \log(n)$, it holds that
    \[
        c_1 \log(n) \leq |\partial_r(\mathcal{S})| \leq 2c_1 \log(n),
    \]
    which implies that there exists an absolute constant $c_2 > 0$ such that 
    \begin{equation}\label{eq-proof-new-prop-2-2}
        |\partial_r(\mathcal{S})| = c_2 \log(n).
    \end{equation}
    Combining \eqref{eq-proof-new-prop-2-1} and \eqref{eq-proof-new-prop-2-2}, we have that
    \[
        \kappa \sqrt{\Delta} /\sigma = c_1c_2^{-1} |\partial_r(
        \mathcal{S})|
    \]
    and therefore $\widetilde{P}_l \in \mathcal{P}$.  The same arguments show that $\widetilde{Q}_l \in \mathcal{P}$, for all $l \in \{1, \ldots, M\}$.
    By construction, we also note that  
    \[
        d_{\mathrm{H}}(\mathcal{S}_{\widetilde{P}_l}, \mathcal{S}_{\widetilde{Q}_m}) \geq n/2, \quad l, m \in \{1, \ldots, M\}.
    \]
    It then follows from Le Cam's lemma \citep[e.g.~Lemma~1 in][]{yu1997assouad} that
    \[
        \inf_{\widehat{\mathcal{S}}} \sup_{P \in \mathcal{P}} \mathbb{E}_P \left\{d_{\mathrm{H}}(\widehat{\mathcal{S}}, \mathcal{S}(P))\right\} \geq \frac{n}{4}\left\{1 - d_{\mathrm{TV}}(\widetilde{P}, \widetilde{Q})\right\},
    \]
    where $d_{\mathrm{TV}}(\cdot, \cdot)$ is the total variation distance between two probability measures and the infimum is over all estimators $\widehat{\mathcal{S}}$.
    
\medskip
\noindent \textbf{Upper bounding the total variation distance}
Let $u_l \in \mathbb{R}^{n/2}$ be a sub-vector of $\widetilde{u}_l$ containing only the first $n/2$ entries.  Let $P_l$ and $P_0$ be the multivariate Gaussian distributions $\mathcal{N}(u_l, \sigma^2 I_{n/2})$ and $\mathcal{N}(0, \sigma^2 I_{n/2})$, respectively.  Due to the symmetry between $\widetilde{u}_l$ and $\widetilde{v}_l$, it holds that
    \[
        d_{\mathrm{TV}}(\widetilde{P}, \widetilde{Q}) \leq 2d_{\mathrm{TV}}(P, P_0), \quad \mbox{where } P = \frac{1}{M} \sum_{m = 1}^M P_l.
    \]
    Since $d_{\mathrm{TV}}(P, P_0) \leq \sqrt{\chi^2(P, P_0))}$ \citep[e.g.~Equation~2.27 in][]{tsybakov2008introduction}, it suffices to provide an upper bound on $\chi^2(P, P_0)$.  We have that
    \begin{align*}
        \chi^2(P, P_0) & = \frac{1}{M^2} \sum_{l, m = 1}^M \mathbb{E}_{P_0}\left(\frac{dP_l dP_m}{dP_0 dP_0}\right) - 1 = \frac{1}{M^2} \sum_{l, m = 1}^M \exp\left(\frac{u_l^{\top} u_m}{\sigma^2}\right) - 1 \\
        & = \frac{1}{M^2} \left[\sum_{l = 1}^M \exp \{c_1 \log(n)\} + M(M-1)\right] -1 = M^{-1} (n^{c_1} - 1).
    \end{align*}
    Recall that $M = n/\{4c_1 \log(n)\}$.  Therefore, for $c_1 \in (0, 1)$, there exists a sufficiently large $n_0 = n(c_1)$ such that for any $n \geq n_0$, $M^{-1} (n^{c_1} - 1) \leq 1/16$.  We therefore complete the proof.
\end{proof}

\begin{proof}[Proof of \Cref{prop-1}]
We prove the result by the Le Cam Lemma \citep[e.g.~Lemma~1 in][]{yu1997assouad}.

\medskip
\noindent \textbf{Constructing distributions}  We consider a complete $n$-node graph $G = (V, E)$, where $V = \{1, \ldots, n\}$ and $E = \{(i, j), 1 \leq i < j \leq n\}$.  

Given an absolute constant $c_1 > 0$, without loss of generality, we assume that $c_1\log(n)$ and $M$ are positive integers, where $M$ satisfies that $2M + 1 = n/\{c_1\log(n)\}$.  For $l \in \{1, \ldots, M\}$, let $\widetilde{u}_l \in \mathbb{R}^n$ be such that the $i$the coordinate of $\widetilde{u}_l(i)$, $i = 1, \ldots, n$, satisfies that
    \[
        \widetilde{u}_l(i) = \begin{cases}
            \sigma, & i \in \{(l-1)c_1 \log(n) + 1, \ldots, l c_1 \log(n)\}, \\
            0, & \mbox{otherwise}.
        \end{cases}
    \]
    Let $\widetilde{v}_l \in \mathbb{R}^n$ be such that $\widetilde{v}_l(i) = \widetilde{u}(n-i+1)$, $i = 1, \ldots, n$.  Let $\widetilde{P}_l$ and $\widetilde{Q}_l$ be multivariate Gaussian distributions $\mathcal{N}(\widetilde{u}_l, \sigma^2 I_n)$ and $\mathcal{N}(\widetilde{v}_l, \sigma^2 I_n)$, respectively and set
    \[
        \widetilde{P} = \frac{1}{M} \sum_{l = 1}^M \widetilde{P}_l \quad \mbox{and} \quad \widetilde{Q} = \frac{1}{M} \sum_{l = 1}^M \widetilde{Q}_l.
    \]

Following the identical arguments as those in the proof of \Cref{prop-2}, we have that $\widetilde{P}_l, \widetilde{Q}_l \in \mathcal{P}$.  By construction, we also note that 
    \[
        d_{\mathrm{H}}(\mathcal{S}_{\widetilde{P}_l}, \mathcal{S}_{\widetilde{Q}_m}) \geq c_1\log(n), \quad l, m \in \{1, \ldots, M\}.
    \]
    It then follows from Le Cam's lemma \citep[e.g.~Lemma~1 in][]{yu1997assouad} that
    \[
        \inf_{\widehat{\mathcal{S}}} \sup_{P \in \mathcal{P}} \mathbb{E}_P \left\{d_{\mathrm{H}}(\widehat{\mathcal{S}}, \mathcal{S}(P))\right\} \geq c_1\log(n)/2 \left\{1 - d_{\mathrm{TV}}(\widetilde{P}, \widetilde{Q})\right\}.
    \]

\medskip
\noindent \textbf{Upper bounding the total variation distance}  Let $u_l \in \mathbb{R}^{n/2}$ be a sub-vector of $\widetilde{u}_l$ containing only the first $n/2$ entries.  Let $P_l$ and $P_0$ be the multivariate Gaussian distributions $\mathcal{N}(u_l, \sigma^2 I_{n/2})$ and $\mathcal{N}(0, \sigma^2 I_{n/2})$, respectively.  It follows from the same arguments as those in the proof of \Cref{prop-2} that
	\begin{align*}
		& d_{\mathrm{TV}}(\widetilde{P}, \widetilde{Q}) \leq 2d_{\mathrm{TV}}(P, P_0) \leq 2\sqrt{\chi^2(P, P_0)} \\
		= & 2\sqrt{\frac{1}{M^2} \sum_{l, m = 1}^M \mathbb{E}_{P_0}\left(\frac{dP_l dP_m}{dP_0 dP_0}\right) - 1} = 2 \sqrt{M^{-1}n^{c_1} - 1},
	\end{align*}
	where the equation follows from the identical arguments as those in the proof of \Cref{prop-2}.  Recall that $M = [n/\{c_1 \log(n)\} - 1]/2$.  Therefore, for any $c_1 \in (0, 1)$, there exists a large enough $n$ such that $d_{\mathrm{TV}}(\widetilde{P}, \widetilde{Q}) < 1/2$ and we conclude the proof.
\end{proof}

\begin{proof}[Proof of \Cref{prop-3}]
Recall that $E = \emptyset$.  We therefore only need to construct nodes in the proof.

\medskip 
\noindent \textbf{$0.2n$-packing number of $2^V$.}  We first construct a collection of distributions.  Let 
	\[
		\mathcal{W} = \left\{A \in 2^V:\, |A| \in [0.49n, 0.51n]\right\}.
	\]
	We introduce independent Rademacher random variables $X_i$ associated with each node $i \in V$.  Therefore
	\begin{align*}
		\mathbb{P}\left\{\left|\sum_{i \in V} X_i - n/2\right| \geq 0.01n \right\} \leq 2\exp (-0.0002n).
	\end{align*}
	This means there exists an absolute constant $c_1 \in [1, 2)$ and a sufficiently large $n_0 \in \mathbb{N}^*$, such that for any $n \geq n_0$,
	\[
		|\mathcal{W}| \geq \{1 - 2\exp (-0.0002n)\} 2^n \geq c_1^n.
	\]

For any $A \in \mathcal{W}$, we let
	\[
		\mathcal{W}_A = \left\{B \in \mathcal{W}: \, |A \triangle B| \leq 0.2n\right\}.
	\]
	Note that for any fixed $A$, 
	\[
		|\mathcal{W}_A| \leq \left|\left\{B \in 2^V: |B| \leq 0.2n \right\}\right|.
	\]
	Using the Rademacher random variables $X_i$'s, we have that
	\[
		\mathbb{P}\left\{\sum_{i \in V} X_i - n/2 \leq -0.3n \right\} \leq \exp (-0.18n).
	\]
	This means for any fixed $A \in \mathcal{W}$ and a sufficiently large $n$, there exists an absolute constant $c_2 \in (0, c_1)$ that
	\[
		|\mathcal{W}_A| \leq \exp (-0.18n) 2^n \leq c_2^n.
	\]
	
We let $\mathcal{M} \subset \mathcal{W}$ and any $A, B \in \mathcal{M}$ satisfy that $|A \triangle B| > 0.2n$.  Then the $0.2n$-packing number of $2^V$ with respect to the Hamming distance satisfies
	\begin{equation}\label{eq-packing-lower}
		M(0.2n, 2^V, d_{\mathrm{Hamm}}) \geq |\mathcal{M}| \geq \frac{|\mathcal{W}|}{|\mathcal{W}_A| + 1} \geq c_3^n,
	\end{equation}
	where $c_3 \in (1, 2)$ is an absolute constant.
	
\medskip
\noindent \textbf{Construction of nodes.}  We now construct a class of distributions $\{P_A, \, A \in \mathcal{M}\}$.  Each $P_A$ is the joint distribution of independent random variables
	\[
		Y_i \sim \begin{cases}
 			\mathcal{N}(\kappa, \sigma^2), & i \in A, \\
 			\mathcal{N}(0, \sigma^2), & i \notin A.
		\end{cases}
	\]	
	We further associate $P_A$ with a connected graph $G_A = (V, E_A)$. 
	
\medskip
\noindent \textbf{Fano's method.}  To use Fano's method, we adopt the version in \cite{yu1997assouad}.  Note that for any $A, B \in \mathcal{M}$, we have that 
	\[
		d_{\mathrm{H}}(A, B) \geq 0.2n,
	\]
	and 
	\[
		\mathrm{KL}(P_A, P_B) = 0.2n \kappa^2/\sigma^2.
	\]
	Then provided that
	\[
		\Delta \kappa^2 \sigma^{-2} \lesssim n,
	\]
	it holds that
	\begin{align*}
		\inf_{\widehat{\mathcal{S}}} \sup_{P \in \mathcal{S}} \mathbb{E}_P \{d_{\mathrm{H}}(\widehat{\mathcal{S}}, \mathcal{S})\} \geq & \frac{0.2n}{2} \left(1 - \frac{0.2n \kappa^2 /\sigma^2 + \log(2)}{\log(M(0.2n, 2^V, d_{\mathrm{Hamm}}))}\right) \\
		\geq & 0.1n \left(1 - \frac{0.2n \kappa^2 /\sigma^2 + \log(2)}{n \log(c_3)}\right) \geq cn.
	\end{align*}
\end{proof}

\begin{proof}[Proof of \Cref{lem-com-resis-weight}]
It follows from Kirchhoff's maximum tree theorem \citep[e.g.~Theorem~1 in][]{chaiken1978matrix} that the number of spanning trees of $G$ equals $\mathrm{det}(L)/n$, where $L$ is the Laplacian of $G$.  Since $G$ is a complete graph, the eigenvalues of $G$'s Laplacian are $n^{n-1}$ and 0. Hence, Kirchhoff's maximum-tree theorem implies that there are $n^{n-2}$ spanning trees.  Due to the definition of spanning trees, each spanning tree consists of $n-1$ edges.  Then there are in total $(n-1)n^{n-2}$ edges contained in all the spanning trees.  

On the other hand, there are $n(n-1)/2$ edges in the complete graph $G$.  Since this is a complete graph, all edges are equivalent.  This implies that each edge $(i, j)$ appears in 
    \[
        \frac{(n-1)n^{n-2}}{n(n-1)/2} = 2 n^{n-3}
    \]
    spanning trees.
    
Due to the definition of the effective resistance weights, we have that
    \[
        \partial_r(\mathcal{S}) = \frac{2 n^{n-3}}{n^{n-2}} n_1 n_2 = 2n_1n_2/n,
    \]
    which concludes the proof.
\end{proof}

\section[second]{PROOFS OF RESULTS IN SECTION~\ref{sec-upper-bound}}

\begin{proof}[Proof of \Cref{thm-1}]
Without loss of generality, in this proof, we assume $f_1^* = 0$, $f_2^* = \kappa$ and $\Delta = |S^*_1| \leq |S^*_2|$.  For any vector $v \in \mathbb{R}^n$ and any subset $A \subset V$, define 
    \[
        \overline{v}_A = |A|^{-1}\sum_{i \in A} v_i.
    \]
    For any partition $\mathcal{C}$, it is associated with a $|\mathcal{C}|$-dimensional subspace $K \subset \mathbb{R}^n$, such that $v \in K$ if and only if $v$ takes a constant value on each element of $\mathcal{C}$.  We denote the orthogonal projection onto $K$ by $P^{\mathcal{C}}: \mathbb{R}^n \to K$.  For the estimator $\{\widehat{S}_1, \widehat{S}_2\}$, we let $A_{kl} = \widehat{S}_k \cap S_l^*$, $k, l = 1, 2$.  Based on this notation, for the uniqueness of the definition, we let
    \begin{equation}\label{eq-def-shat-1-2}
        \frac{|A_{11}|}{|\widehat{S}_1|} \geq \frac{|A_{21}|}{|\widehat{S}_2|}.
    \end{equation}
            
\medskip
\noindent \textbf{Step 1.}  Let 
	\begin{equation}\label{eq-large-prob-event-fan}
	    \mathcal{E} = \left\{\|P^{\mathcal{C}} \varepsilon\|^2 \leq C |\partial_w(\mathcal{C})| \log\{w(E)\}, \, \forall \mbox{ partition } \mathcal{C} \mbox{ of } G\right\}.
	\end{equation}
	It follows from Lemma~B.2 in \cite{fan2018approximate} that
	\[
		\mathbb{P}\{\mathcal{E}\} \geq 1 - w(E)^{-c},
	\]
	where $C, c > 0$ are absolute constants.  The rest of the proof is conducted on the event $\mathcal{E}$. 

\medskip
\noindent \textbf{Step 2.}  Let $\{\widehat{S}_1, \widehat{S}_2\}$ be the output of \Cref{alg-main}.   It must hold that
	\begin{align*}
		& \sum_{i \in \widehat{S}_1} (Y_i - \overline{Y}_{\widehat{S}_1}^{\delta})^2 + \sum_{i \in \widehat{S}_2} (Y_i - \overline{Y})^2 + 2\lambda |\partial_w (\widehat{S}_1, \widehat{S}_2)| \\
		\leq & \sum_{i \in S_1^*} (Y_i - \overline{Y}_1^{\delta})^2 + \sum_{i \in S_2^*} (Y_i - \overline{Y})^2 + 2\lambda |\partial_w (S_1^*, S_2^*)| + 4\tau,
	\end{align*}
	which implies that
	\begin{align}
		& \sum_{i \in \widehat{S}_1} (Y_i - \overline{Y}_{\widehat{S}_1})^2 + \sum_{i \in \widehat{S}_2} (Y_i - \overline{Y})^2 + 2\lambda |\partial_w (\widehat{S}_1, \widehat{S}_2)|  \nonumber \\
		\leq & \sum_{i \in S_1^*} (Y_i - \overline{Y}_1)^2 + \sum_{i \in S_2^*} (Y_i - \overline{Y})^2 + 2\lambda |\partial_w (S_1^*, S_2^*)| + 2n\delta^2 + 4\tau. \label{eq-main-eq-in-the-proof}
	\end{align}
	
Let 
	\[
		(\mu_1)_i = \begin{cases}
 			\overline{Y}_{\widehat{S}_1}, & i \in \widehat{S}_1, \\
 			\overline{Y}, & i \in \widehat{S}_2,
		\end{cases} \quad
		(f_1)_i = \begin{cases}
 			\frac{|A_{12}|}{|\widehat{S}_1|}\kappa, & i \in \widehat{S}_1,\\
 			\frac{|S_2^*|}{n} \kappa, & i \in \widehat{S}_2,
 		\end{cases} \quad
 		(\varepsilon_1)_i = \begin{cases}
 		    \frac{1}{|\widehat{S}_1|} \sum_{i \in \widehat{S}_1} \varepsilon_i, & i \in \widehat{S}_1, \\
 		    \frac{1}{n} \sum_{i \in V} \varepsilon_i, & i \in \widehat{S}_2,
 		\end{cases}
	\]
	\[
		(\mu_2)_i = \begin{cases}
 			\overline{Y}_{S_1^*}, & i \in S_1^*,\\
 			\overline{Y}, & i \in S_2^*,
 		\end{cases} \quad
 		(f_2)_i = \begin{cases}
 			0, & i \in S_1^*, \\
 			\frac{|S_2^*|}{n} \kappa, & i \in S_2^*,  			
 		\end{cases} \quad
 		\varepsilon_2 = \begin{cases}
 		    \frac{1}{|S^*_1|} \sum_{i \in S^*_1} \varepsilon_i, & i \in S^*_1, \\
 		    \frac{1}{n} \sum_{i \in V} \varepsilon_i, & i \in S^*_2, 
 		\end{cases}
	\]
	\[
		(\mu_0)_i = \begin{cases}
 			0, & i \in S_1^*, \\
 			\kappa, & i \in S_2^*
 		\end{cases} \quad \mbox{and} \quad Y = \mu_0 + \varepsilon,
	\]
	where, we have assumed that $\widehat{S}_1 \neq \emptyset$.  We will come back to prove this claim in \textbf{Step 5}.
	
Based on the notation above, we note that 	
	\begin{align*}
		(f_1 + f_2 - 2\mu_0)_i = \begin{cases}
 			\frac{|A_{12}|}{|\widehat{S}_1|} \kappa, & i \in A_{11},\\
	 		- \left(\frac{|A_{11}|}{|\widehat{S}_1|} + \frac{|S_1^*|}{n}\right) \kappa, & i \in A_{12},\\
 			\frac{|S_2^*|}{n}\kappa, & i \in A_{21}, \\
 			-2\frac{|S_1^*|}{n}\kappa, & i \in A_{22},
	 	\end{cases} \quad
 		(f_1 - f_2)_i = \begin{cases}
 			\frac{|A_{12}|}{|\widehat{S}_1|}\kappa, & i \in A_{11}, \\
	 		\left(\frac{|A_{12}|}{|\widehat{S}_1|} - \frac{|S_2^*|}{n}\right) \kappa, & i \in A_{12}, \\
 			\frac{|S_2^*|}{n} \kappa, & i \in A_{21}, \\
 			0, & i \in A_{22}.
	 	\end{cases}
	\end{align*}
	and
	\[
		(\varepsilon_1 - \varepsilon_2)_i = \begin{cases}
 			\frac{1}{|\widehat{S}_1|}\sum_{j \in \widehat{S}_1} \varepsilon_j - \frac{1}{|S_1^*|} \sum_{j \in S_1^*} \varepsilon_j, & i \in A_{11}, \\
	 		\frac{|\widehat{S}_2|}{n|\widehat{S}_1|}\sum_{j \in \widehat{S}_1} \varepsilon_j - \frac{1}{n} \sum_{j \in \widehat{S}_2} \varepsilon_j , & i \in A_{12}, \\
	 		\frac{1}{n} \sum_{j \in S^*_2} \varepsilon_j - \frac{|S^*_2|}{n|S^*_1|} \sum_{j \in S^*_1} \varepsilon_j, & i \in A_{21}, \\
	 		0, & i \in A_{22}.
	 	\end{cases}
	\]
	
With the above notation, we have that $\mu_1 = f_1 + \varepsilon_1$, $\mu_2 = f_2 + \varepsilon_2$,
    \[
        \sum_{i \in \widehat{S}_1} (Y_i - \overline{Y}_{\widehat{S}_1})^2 + \sum_{i \in \widehat{S}_2} (Y_i - \overline{Y})^2 = \|Y - \mu_1\|^2 \quad \mbox{and} \quad \sum_{i \in S_1^*} (Y_i - \overline{Y}_1)^2 + \sum_{i \in S_2^*} (Y_i - \overline{Y})^2 = \|Y - \mu_2\|^2.
    \]
    
\medskip
\noindent \textbf{Step 3.}  In this step, we are to exploit \eqref{eq-def-shat-1-2}, which directly implies that
    \[
        |A_{11}| |A_{22}| > |A_{12}||A_{21}|.
    \]
    This means 
    \[
        |A_{11}| |A_{22}| + |A_{11}||A_{12}| + |A_{12}|^2 + |A_{12}||A_{22}| > |A_{12}||A_{21}| + |A_{11}||A_{12}| + |A_{12}|^2 + |A_{12}||A_{22}|,
    \]
    which is equivalent to
    \[
        (|A_{11}| + |A_{12}|)(|A_{12}| + |A_{22}|) > |A_{12}|(|A_{11}| + |A_{12}| + |A_{21}| + |A_{22}|).
    \]
    Simplifying the above leads to
    \begin{equation}\label{eq-new-proof-key}
        \frac{|S^*_2|}{n} > \frac{|A_{12}|}{|\widehat{S}_1|}.
    \end{equation}
    
Another consequence of \eqref{eq-def-shat-1-2} is that 
    \begin{equation}\label{eq-new-proof-key-22}
        \frac{|A_{11}|}{|\widehat{S}_1|} \geq \frac{|S^*_1|}{n}.
    \end{equation}
    If \eqref{eq-new-proof-key-22} does not hold, then we have that
    \[
        \frac{|S^*_1|}{n} = \frac{|A_{11}| + |A_{21}|}{|\widehat{S}_1| + |\widehat{S}_2|} \leq \frac{|A_{11}|}{|\widehat{S}_1|} < \frac{|S^*_1|}{n},
    \]
    which is a contradiction

\medskip  
\noindent \textbf{Step 4.}	In this step, we are to lower bound
	\[
		Q = \|Y - \mu_1\|^2 - \|Y - \mu_2\|^2.
	\]	
	
Note that, with an absolute constant $c \in (0, 1)$, it holds that 
	\begin{align}\label{eq-q-0}
		Q & = \|\mu_0 - \mu_1 + \varepsilon\|^2 - \|\mu_0 - \mu_2 + \varepsilon\|^2 = \|\mu_0 - \mu_1\|^2 - \|\mu_0 - \mu_2\|^2 + 2(\mu_2 - \mu_1)^{\top}\varepsilon \nonumber \\
		& \geq  \|\mu_0 - f_1 - \varepsilon_1\|^2 - 	\|\mu_0 - f_2 - \varepsilon_2\|^2 - c\|f_1 - f_2 + \varepsilon_1 - \varepsilon_2\|^2 - \frac{1}{c} \|P^{\widehat{\mathcal{S}} \vee \mathcal{S}^*}\varepsilon\|^2 \nonumber \\
		& = (f_1 + f_2 - 2\mu_0)^{\top}(f_1 - f_2) + \langle f_1 + f_2 - 2\mu_0, \varepsilon_1 - \varepsilon_2\rangle + \langle f_1 - f_2, \varepsilon_1 + \varepsilon_2\rangle + \langle \varepsilon_1 \nonumber \\
		& \hspace{1cm} + \varepsilon_2, \varepsilon_1 - \varepsilon_2\rangle - c\|f_1 - f_2 + \varepsilon_1 - \varepsilon_2\|^2 - \frac{1}{c} \|P^{\widehat{\mathcal{S}} \vee \mathcal{S}^*}\varepsilon\|^2 \nonumber \\
		& \geq (f_1 + f_2 - 2\mu_0)^{\top}(f_1 - f_2) +  \langle f_1 + f_2 - 2\mu_0, \varepsilon_1 - \varepsilon_2\rangle - c_1\|f_1 - f_2\|^2 - \frac{1}{4c_1} \|P^{\widehat{\mathcal{S}} \vee \mathcal{S}^*} (\varepsilon_1 + \varepsilon_2)\|^2 \nonumber \\
		& \hspace{1cm} - c_1\|\varepsilon_1 - \varepsilon_2\|^2 - \frac{1}{4c_1} \|P^{\widehat{\mathcal{S}} \vee \mathcal{S}^*} (\varepsilon_1 + \varepsilon_2)\|^2 - 2c \|f_1 - f_2\|^2 - 2c\|\varepsilon_1 - \varepsilon_2\|^2 - \frac{1}{c}\|P^{\widehat{\mathcal{S}} \vee \mathcal{S}^*} \varepsilon\|^2 \nonumber \\
		& = (f_1 + f_2 - 2\mu_0)^{\top}(f_1 - f_2) +  \langle f_1 + f_2 - 2\mu_0, \varepsilon_1 - \varepsilon_2\rangle - (c_1 + 2c)\|f_1 - f_2\|^2 \nonumber \\
		& \hspace{1cm} - \left(\frac{2}{c_1} + \frac{1}{c}\right) \|P^{\widehat{\mathcal{S}} \vee \mathcal{S}^*} \varepsilon\|^2 - (c_1 + 2c) \|\varepsilon_1 - \varepsilon_2\|^2,
	\end{align}
	where facts that $\mu_1 = f_1 + \varepsilon_1$, $\mu_2 = f_2 + \varepsilon_2$ and $2ab \geq ca^2 + b^2/c$, with $c > 0$, are repeatedly used above.

For the first term in \eqref{eq-q-0}, we have that 
	\begin{align}
		& (f_1 + f_2 - 2\mu_0)^{\top}(f_1 - f_2) \nonumber \\
		= & \frac{|A_{11}||A_{12}|^2}{|\widehat{S}_1|^2} \kappa^2 - \left(\frac{|A_{11}||A_{12}|^2}{|\widehat{S}_1|^2} \kappa^2 - \frac{|A_{11}||A_{12}||S_2^*|}{|\widehat{S}_1| n} \kappa^2 + \frac{|A_{12}|^2|S_1^*|}{|\widehat{S}_1|n}\kappa^2 - \frac{|A_{12}||S_1^*||S_2^*|}{n^2} \kappa^2 \right)\label{eq-q-1-11111}\\
		& \hspace{1cm} + \frac{|A_{21}||S_2^*|^2}{n^2} \kappa^2 \label{eq-q-1-21111} \\
		= & \frac{|A_{11}||A_{12}||S_2^*|}{|\widehat{S}_1| n} \kappa^2 - \frac{|A_{12}|^2|S_1^*|}{|\widehat{S}_1|n}\kappa^2 + \frac{|A_{12}||S_1^*||S_2^*|}{n^2} \kappa^2 + \frac{|A_{21}||S_2^*|^2}{n^2} \kappa^2 \nonumber  \\
		\geq & \frac{|A_{11}|}{2|\widehat{S}_1|} |A_{12}|\kappa^2 + \frac{|A_{12}||S^*_1|}{n}\left(\frac{|S^*_2|}{n} - \frac{|A_{12}|}{|\widehat{S}_1|}\right) \kappa^2 + \frac{|A_{21}|}{4} \kappa^2 \nonumber \\
		\geq & \left(\frac{|A_{11}|}{2|\widehat{S}_1|} |A_{12}| + \frac{1}{4}|A_{21}|\right)\kappa^2, \label{eq-new-proof-12121}
	\end{align}
	where the first term in \eqref{eq-q-1-11111} is from the products indexed in $A_{11}$, the four terms in the brackets in \eqref{eq-q-1-11111} are from the products indexed in $A_{12}$, the term in \eqref{eq-q-1-21111} is from the products indexed in $A_{21}$, the third inequality is due to the fact that $|S^*_2| \geq |S^*_1|$ and the final inequality is due to \eqref{eq-new-proof-key}.

For the second term in \eqref{eq-q-0}, with the notation that 
    \[
		\overline{\varepsilon}_{kl} = |A_{kl}|^{-1}\sum_{i \in A_{kl}} \varepsilon_i, \quad k, l = 1, 2,
	\]
    we have that
    \begin{align}
		& \langle f_1 + f_2 - 2\mu_0, \varepsilon_1 - \varepsilon_2\rangle \nonumber \\
		= & \frac{|A_{11}||A_{12}|\kappa}{|\widehat{S}_1|} \left(\frac{1}{|\widehat{S}_1|}\sum_{j \in \widehat{S}_1} \varepsilon_i - \frac{1}{|S_1^*|} \sum_{j \in S_1^*} \varepsilon_i\right) - |A_{12}|\left(\frac{|A_{11}|}{|\widehat{S}_1|} + \frac{|S_1^*|}{n}\right) \kappa \left(\frac{|\widehat{S}_2|}{n|\widehat{S}_1|}\sum_{j \in \widehat{S}_1} \varepsilon_j - \frac{1}{n} \sum_{j \in \widehat{S}_2} \varepsilon_j \right)  \nonumber \\
		& + \frac{|S_2^*||A_{21}|}{n}\kappa \left(	\frac{1}{n} \sum_{j \in S^*_2} \varepsilon_j - \frac{|S^*_2|}{n|S^*_1|} \sum_{j \in S^*_1} \varepsilon_j\right)  \nonumber \\
		\leq & \left|\frac{|A_{11}|^2 |A_{12}|^2}{|\widehat{S}_1|^2 |S^*_1|} - \frac{|A_{11}|^2 |A_{21}| |A_{12}|}{|\widehat{S}_1|^2 |S^*_1|} - \frac{|A_{12}| |A_{11}|^2 |\widehat{S}_2|}{n |\widehat{S}_1|^2} - \frac{|A_{12}| |S_1^*| |\widehat{S}_2| |A_{11}|}{n^2 |\widehat{S}_1|} - \frac{|S^*_2|^2 |A_{21}| |A_{11}|}{n^2 |S_1^*|}\right| \kappa |\overline{\varepsilon}_{11}|  \nonumber \\
		& + \left|\frac{|A_{12}| |A_{11}| |A_{22}|}{n |\widehat{S}_1|} + \frac{|A_{12}| |S^*_1| |A_{22}|}{n^2} + \frac{|S^*_2||A_{21}||A_{22}|}{n^2}\right| \kappa |\overline{\varepsilon}_{22} |  \nonumber \\
		& + \left|\frac{|A_{11}| |A_{12}|^2}{|\widehat{S}_1|^2} - \frac{|A_{12}|^2 |A_{11}| |\widehat{S}_2|}{n |\widehat{S}_1|^2} - \frac{|A_{12}|^2 |S^*_1| |\widehat{S}_2|}{n^2 |\widehat{S}_1|} + \frac{|S^*_2| |A_{21}| |A_{12}|}{n^2}\right|\kappa |\overline{\varepsilon}_{12}| \nonumber \\
		& + \left|-\frac{|A_{11}||A_{12}||A_{21}|}{|\widehat{S}_1||S^*_1|} + \frac{|A_{12}| |A_{11}| |A_{21}|}{n |\widehat{S}_1|} + \frac{|A_{12}| |A_{21}| |S^*_1|}{n^2} - \frac{|S^*_2|^2 |A_{21}|^2}{n^2 |S^*_1|}\right|\kappa |\overline{\varepsilon}_{21}| \nonumber \\
		\leq & c B \kappa^2 + C \|P^{\widehat{\mathcal{S}} \vee \mathcal{S}^*} \varepsilon\|^2, \label{eq-new-proof-1211111}
	\end{align}
	where $c, C > 0$ are absolute constants and
	\begin{align}
	    B & = \frac{|A_{11}|^4 |A_{12}|^4}{|\widehat{S}_1|^4 |S^*_1|^2 |A_{11}|} + \frac{|A_{11}|^4 |A_{21}|^2 |A_{12}|^2}{|\widehat{S}_1|^4 |S^*_1|^2 |A_{11}|} + \frac{|A_{12}|^2 |A_{11}|^4 |\widehat{S}_2|^2}{n^2 |\widehat{S}_1|^4 |A_{11}|} + \frac{|A_{12}|^2 |S_1^*|^2 |\widehat{S}_2|^2 |A_{11}|^2}{n^4 |\widehat{S}_1|^2|A_{11}|} + \frac{|S^*_2|^4 |A_{21}|^2 |A_{11}|^2}{n^4 |S_1^*|^2 |A_{11}|} \nonumber \\
	    & \hspace{0.3cm} + \frac{|A_{12}|^2 |A_{11}|^2 |A_{22}|^2}{n^2 |\widehat{S}_1|^2 |A_{22}|} + \frac{|A_{12}|^2 |S^*_1|^2 |A_{22}|^2}{n^4 |A_{22}|} + \frac{|S^*_2|^2|A_{21}|^2|A_{22}|^2}{n^4 |A_{22}|} \nonumber \\
	    & \hspace{0.3cm} + \frac{|A_{11}|^2 |A_{12}|^4}{|\widehat{S}_1|^4 |A_{12}|} + \frac{|A_{12}|^4 |A_{11}|^2 |\widehat{S}_2|^2}{n^2 |\widehat{S}_1|^4 |A_{12}|} + \frac{|A_{12}|^4 |S^*_1|^2 |\widehat{S}_2|^2}{n^4 |\widehat{S}_1|^2 |A_{12}|} + \frac{|S^*_2|^2 |A_{21}|^2 |A_{12}|^2}{n^4 |A_{12}|} \nonumber \\
	    & \hspace{0.3cm} + \frac{|A_{11}|^2|A_{12}|^2|A_{21}|^2}{|\widehat{S}_1|^2|S^*_1|^2 |A_{21}|} + \frac{|A_{12}|^2 |A_{11}|^2 |A_{21}|^2}{n^2 |\widehat{S}_1|^2 |A_{21}|} + \frac{|A_{12}|^2 |A_{21}|^2 |S^*_1|^2}{n^4 |A_{21}|} + \frac{|S^*_2|^4 |A_{21}|^4}{n^4 |S^*_1|^2 |A_{21}|} \nonumber \\
	    & \leq \frac{|A_{11}|}{|\widehat{S}_1|}|A_{12}| + \frac{|A_{11}|^3}{|\widehat{S}_1|^3}|A_{12}| + \frac{|A_{11}|^3}{|\widehat{S}_1|^3}|A_{12}| + \frac{|A_{11}|}{|\widehat{S}_1|}|A_{12}| + |A_{21}| \nonumber \\
	    & \hspace{0.3cm} + \frac{|A_{11}|^2}{|\widehat{S}_1|^2}|A_{12}| + \frac{|S^*_1|}{n} |A_{12}| + |A_{21}| \nonumber \\
	    & \hspace{0.3cm} + \frac{|A_{11}|^3}{|\widehat{S}_1|^3}|A_{12}| + \frac{|A_{11}|^2}{|\widehat{S}_1|^2}|A_{12}| + \frac{|S^*_1|^2}{n^2} |A_{12}| + |A_{21}| \nonumber \\
	    & \hspace{0.3cm} + |A_{21}| +  \frac{|A_{11}|}{|\widehat{S}_1|}|A_{12}| + \frac{|S^*_1|^2}{n^2} |A_{12}| + |A_{21}| \nonumber \\
	    & \leq \frac{11|A_{11}|}{|\widehat{S}_1|} |A_{12}| + 5 |A_{21}|,  \label{eq-new-proof-999}
    \end{align}
    where the last inequality is due to \eqref{eq-new-proof-key-22}.

For the third term in \eqref{eq-q-0}, we have that 
	\begin{align} \label{eq-q-2}
		& \|f_1 - f_2\|^2 = \frac{|A_{12}|^2 |A_{11}|}{|\widehat{S}_1|^2}\kappa^2 +  \left(\frac{|A_{12}|}{|\widehat{S}_1|} - \frac{|S_2^*|}{n}\right)^2 |A_{12}| \kappa^2 + \frac{|S^*_2|^2|A_{21}|}{n^2}\kappa^2 \nonumber \\
		\leq & \frac{|A_{11}|}{|\widehat{S}_1|} |A_{12}|\kappa^2 + \left(\frac{|S^*_1|}{n} - \frac{|A_{11}|}{|\widehat{S}_1|}\right)^2  |A_{12}| \kappa^2 + \frac{1}{4}|A_{21}| \kappa^2 \nonumber \\
		\leq & \frac{5|A_{11}|}{|\widehat{S}_1|} |A_{12}|\kappa^2 + \frac{1}{4}|A_{21}| \kappa^2,
	\end{align}
    where the last inequality is due to \eqref{eq-new-proof-key-22}.

For the last term in \eqref{eq-q-0}, we have that 
    \begin{align}
        & \|\varepsilon_1 - \varepsilon_2\|^2 \nonumber \\
        \leq & |A_{11}|\left(\frac{|A_{12}||A_{11}| - |A_{21}||A_{11}|}{|\widehat{S}_1| |S^*_1|} \overline{\varepsilon}_{11} + \frac{|A_{12}|}{|\widehat{S}_1|}\overline{\varepsilon}_{12} - \frac{|A_{21}|}{|S^*_1|}\overline{\varepsilon}_{21}\right)^2 \nonumber \\
        & \hspace{1cm} + |A_{12}| \left(\frac{|\widehat{S}_2||A_{11}|}{n|\widehat{S}_1|} \overline{\varepsilon}_{11} + \frac{|\widehat{S}_2||A_{12}|}{n|\widehat{S}_1|} \overline{\varepsilon}_{12} - \frac{|A_{21}|}{n} \overline{\varepsilon}_{21} - \frac{|A_{22}|}{n} \overline{\varepsilon}_{22}\right)^2 \nonumber \\
        & \hspace{1cm} + |A_{21}| \left(-\frac{|S^*_2||A_{11}|}{n|S^*_1|} \overline{\varepsilon}_{11} + \frac{|A_{12}|}{n}\overline{\varepsilon}_{12} - \frac{|S^*_2||A_{21}|}{n|S^*_1|}\overline{\varepsilon}_{21} + \frac{|A_{22}|}{n} \overline{\varepsilon}_{22}\right)^2  \nonumber \\
        \leq & 4\left(\frac{|A_{11}|^3|A_{12}|^2}{|\widehat{S}_1|^2 |S^*_1|^2} + \frac{|A_{11}|^3 |A_{21}|^2}{|\widehat{S}_1|^2 |S^*_1|^2} + \frac{|A_{12}||\widehat{S}_2|^2|A_{11}|^2}{n^2|\widehat{S}_1|^2} + \frac{|A_{21}| |S^*_2|^2|A_{11}|^2}{n^2 |S^*_1|^2}\right) \overline{\varepsilon}_{11}^2 \nonumber \\
        & \hspace{1cm} + 4 \left(\frac{|A_{11}||A_{12}|^2}{|\widehat{S}_1|^2} + \frac{|A_{12}| |\widehat{S}_2|^2 |A_{12}|^2}{n^2 |\widehat{S}_1|^2} + \frac{|A_{21}| |A_{12}|^2}{n^2} \right)\overline{\varepsilon}_{12}^2 \nonumber \\
        & \hspace{1cm} + 4 \left(\frac{|A_{11}||A_{21}|^2}{|S^*_1|^2} + \frac{|A_{12}| |A_{21}|^2}{n^2} + \frac{|A_{21}||S^*_2|^2|A_{21}|^2}{n^2|S^*_1|^2}\right) \overline{\varepsilon}_{21}^2 + 4 \left(\frac{|A_{12}||A_{22}|^2}{n^2} + \frac{|A_{21}||A_{22}|^2}{n^2}\right) \overline{\varepsilon}_{22}^2 \nonumber \\
        \leq & 16 |A_{11}|\overline{\varepsilon}_{11}^2 + 12|A_{12}| \overline{\varepsilon}_{12}^2 + 12|A_{21}| \overline{\varepsilon}_{21}^2 + 8|A_{22}| \overline{\varepsilon}_{22}^2 \nonumber \\
        \leq & 16 \|P^{\widehat{\mathcal{S}} \vee \mathcal{S}^*} \varepsilon\|^2. \label{eq-aaaaaaaaaa}
    \end{align}

Combining \eqref{eq-q-0}, \eqref{eq-new-proof-12121}, \eqref{eq-new-proof-1211111}, \eqref{eq-new-proof-999}, \eqref{eq-q-2} and \eqref{eq-aaaaaaaaaa}, we have that, with absolute constants $c_*, C_1 > 0$, 
    \begin{align}\label{eq-q-lb}
        Q & \geq c_* \left(\frac{|A_{11}| |A_{12}|}{|\widehat{S}_1|} + |A_{21}|\right) \kappa^2 - C_1\sigma^2 |\partial_w(S_1^*, S_2^*)| - C_1\sigma^2|\partial_w(\widehat{S}_1, \widehat{S}_2)|,
	\end{align}
	where the second inequality is due to the choices of the constants.

Combining \eqref{eq-main-eq-in-the-proof} and \eqref{eq-q-lb}, with a sufficiently large $C_{\lambda} > 0$, we have that 
	\[
		C\sigma^2|\partial_w(\widehat{S}_1, \widehat{S}_2)| + c_* \left(\frac{|A_{11}| |A_{12}|}{|\widehat{S}_1|} + |A_{21}|\right) \kappa^2 \leq 2\lambda|\partial_w(S_1^*, S^*_2)| + 2n\delta^2 + \tau,
	\]
	which implies \eqref{eq-thm1-result}.

\medskip
\noindent\textbf{Step 5.} To conclude the proof, the only remaining task is to show that $\widehat{S}_1 \neq \emptyset$.  We prove by contradiction and therefore have that
	\[
		\Delta (\overline{Y} - \overline{Y}_1)^2 \leq 2\lambda|\partial_w(\mathcal{S}^*)| + 2n\delta^2+ \tau.
	\]
	On the event $\mathcal{E}$, it implies that 
	\begin{align*}
		\Delta \kappa^2 \leq C\lambda |\partial_w (\mathcal{S})| + 2n \delta^2+ \tau,
	\end{align*}
	which contradicts with \Cref{assume-snr}.  We thus conclude the proof.
\end{proof}

\begin{proof}[Proof of \Cref{prop-constant}]
Note that due to the design of the algorithms, the output of \Cref{alg-main} has fewer pieces than the output of Algorithm~1 in \cite{fan2018approximate}.  Therefore, the claim holds directly follows from Theorem~3.5 in \cite{fan2018approximate}.
\end{proof}

\end{document}